\documentclass[oneside,13pt]{article}
\usepackage{amsfonts}

\usepackage{jmlr2e}
\usepackage{hyperref}
\usepackage{amsmath}
\usepackage{algorithm}
\usepackage{algorithmic}

\usepackage{color}
\newtheorem{assumption}{Assumption}

\usepackage{graphicx}
\usepackage[tight,footnotesize]{subfigure}
\begin{document}

\title{Stochastic Alternating Direction Method of Multipliers with Variance Reduction for Nonconvex Optimizations }

\author{\name Feihu Huang \email huangfeihu@nuaa.edu.cn \\
       \addr College of Computer Science and Technology\\
       Nanjing University of Aeronautics and Astronautics, Nanjing, 210016, China
       \AND
       \name Songcan Chen \email s.chen@nuaa.edu.cn \\
       \addr College of Computer Science and Technology\\
       Nanjing University of Aeronautics and Astronautics, Nanjing, 210016, China
       \AND
       \name Zhaosong Lu \email zhaosong@sfu.ca \\
       \addr Department of Mathematics, \\
       Simon Fraser University, Burnaby, BC V5A 1S6, Canada}

\editor{ {\color{blue}{ The first version is provided in October 9, 2016. }} }

\maketitle
\begin{abstract}
In the paper, we study the stochastic alternating direction method of multipliers (ADMM)
for the nonconvex optimizations, and propose three classes of the nonconvex stochastic ADMM with variance reduction,
based on different reduced variance stochastic gradients.
Specifically, the first class called the nonconvex stochastic variance reduced gradient ADMM (SVRG-ADMM),
uses a multi-stage scheme
to progressively reduce the variance of stochastic gradients.
The second is the nonconvex stochastic average gradient ADMM (SAG-ADMM),
which additionally uses the old gradients estimated in the previous iteration.
The third called SAGA-ADMM is an extension of the SAG-ADMM method.
Moreover, under some mild conditions, we establish the iteration complexity bound of $O(1/\epsilon)$ of the proposed methods
to obtain an $\epsilon$-stationary solution of the nonconvex optimizations.
In particular, we provide a general framework to analyze the iteration complexity
of these nonconvex stochastic ADMM methods with variance reduction.
Finally, some numerical experiments demonstrate the effectiveness of our methods.

\end{abstract}

\section{Introduction}
Stochastic optimization method is a class of powerful optimization tool for solving large-scale problems in machine learning.
For example, the stochastic gradient descent (SGD) \citep{bottou2004stochastic} is an efficient method
for solving the finite-sum optimization problem,
which is a fundamental to machine learning.
Specifically, the SGD only computes gradient of one sample instead of visiting all samples in each iteration.
Though its scalability, due to the variance in the stochastic process, the SGD has slower convergence rate
than the batch gradient method.
Recently, many accelerated versions of the SGD have successfully been proposed to reduce the variances, 
and obtain some better convergence rates.
For example, the stochastic average gradient (SAG) method \citep{roux2012stochastic} obtains a fast convergence rate
by incorporating the old gradients estimated in the previous iterations.
The stochastic dual coordinate ascent (SDCA) method \citep{shalev2013stochastic} performs the stochastic coordinate ascent on the dual problems
and also obtains a fast convergence rate.
Moreover, an accelerated randomized proximal coordinate gradient method (APCG) \citep{lin2015accelerated}
accelerates the SDCA method by using the Nesterov's acceleration technique \citep{nesterov2004introductory}.
However, these accelerated methods obtain faster convergence rate than the standard SGD
at the cost of requiring much space to store old gradients or dual variables.
To deal with this dilemma, thus, the stochastic variance reduced gradient (SVRG) methods \citep{johnson2013accelerating,xiao2014proximal} are proposed,
and enjoy a fast convergence rate with no extra space to store the intermediate gradients or dual variables.
Moreover, \cite{defazio2014saga} have proposed a novel method called SAGA, 
which extends the SAG method and enjoys better theoretical convergence rate
than both the SAG and SVRG methods.

Though the above gradient-based methods can effectively solve many problems in machine learning,
they are still difficultly competent for some complicated problems,
such as the graph-guided SVM \citep{ouyang2013stochastic} and the latent variable graphical models \citep{ma2013alternating}.
It is well known that the alternating direction method
of multipliers (ADMM) \citep{gabay1976dual,boyd2011distributed} has been advocated
as an efficient optimization method in many application fields
such as machine learning \citep{danaher2014joint} and statistics \citep{fang2015generalized}.
However, the offline or batch  ADMM need to compute an empirical risk loss function on all training samples
at each iteration, which makes it unsuitable for large-scale learning problems.
Thus, the online or stochastic versions of ADMM \citep{wang2012online,suzuki2013dual,ouyang2013stochastic} have been developed
for the large-scale/stochastic optimizations.
Due to the variance in the stochastic process, these initial stochastic ADMM methods also suffer from slow convergence rate.
Recently, some accelerated stochastic ADMM methods are proposed to efficiently solve the large-scale learning problems.
For example, a fast stochastic ADMM \citep{zhong2014fast} is proposed via incorporating the previous estimated gradients.
\cite{azadi2014towards} has proposed an accelerated stochastic ADMM
by using Nesterov's accelerated method \citep{nesterov2004introductory}.
Moreover, an adaptive stochastic ADMM \citep{zhao2015adaptive} is proposed by using the adaptive stochastic gradients.
The stochastic dual coordinate ascent ADMM \citep{suzuki2014stochastic} obtains a fast convergence rate by solving the dual problems.
More recently, the scalable stochastic ADMMs \citep{zhao2015scalable,zheng2016fast} are developed,
and obtain fast convergence rates with no extra space for the previous gradients or dual variables.

So far, the above study on stochastic optimization methods relies heavily on the strongly convex or convex objective functions.
However, there exist many useful nonconvex models in machine learning such as
the nonconvex robust empirical risk minimization models \citep{aravkin2016smart}
and deep learning \citep{lecun2015deep}.
Thus, the study of stochastic optimization methods for the nonconvex problems is much needed.
Recently, some works focus on studying the stochastic gradient methods for the nonconvex optimizations.
For example, \cite{ghadimi2016accelerated} and \cite{ghadimi2016mini} have established the iteration complexity of $O(1/\epsilon^2)$ for the SGD
to obtain an $\epsilon$-stationary solution of the nonconvex optimizations.
\cite{allen2016variance,reddi2016stochastic,reddi2016fast1} have proved that both the nonconvex SVRG and SAGA methods
obtain an iteration complexity of $O(1/\epsilon)$ for the nonconvex optimizations.
In particular, \cite{li2016stochastic} have studied the stochastic gradient method for nonconvex sparse learning via variance reduction,
which reaches an asymptotically linear convergence rate by exploring the properties of the specific problems.
Moreover, \cite{reddi2016fast2,aravkin2016smart} have studied the variance reduced stochastic methods
for the nonconvex nonsmooth composite problems,
and have proved that they have the iteration complexity of $O(1/\epsilon)$.
At the same time, \cite{hajinezhad2016nestt} have proposed a nonconvex distributed and stochastic primal dual splitting method
for the nonconvex nonsmooth problems
and prove that it also has the iteration complexity of $O(1/\epsilon)$ to obtain an $\epsilon$-stationary solution.

Similarly, the above nonconvex methods are difficult to be competent
to some complicated nonconvex problems, such as the  graph-guided regularization risk loss minimizations
and tensor decomposition \citep{jiang2016structured}.
Fortunately, it is found that the ADMM method can be well competent to these complicated nonconvex problems,
though it may fail to converge sometimes.
Recently, some work \citep{wang2015convergence,yang2015alternating,wang2015global,hong2016convergence,jiang2016structured}
begin to devote to the study of the ADMM method for the nonconvex optimizations.
However, they mainly focus on studying the determinate ADMM for the nonconvex optimizations.
Due to computing the empirical loss function on
all the training examples at each iteration, these nonconvex ADMMs can not be well competent to the large-scale learning problems.
Though \cite{hong2014distributed} has proposed a distributed, asynchronous and incremental algorithm based on the ADMM method
for the large-scale nonconvex problems, the proposed method is still difficult to be competent
to these complicated nonconvex problems such as the graph-guided models, and
its iteration complexity is not provided.
At present, to the best of our knowledge, there still exists few study of the stochastic ADMM for the noncovex optimizations.
In the paper, thus, we study the stochastic ADMM methods for solving the nonconvex nonsmooth stochastic optimizations as follows:
\begin{align}
 & \min_{x,y} \mathbb{E}_{\xi}[F(x,\xi)]+g(y)  \\
 & \mbox{s.t.} \ Ax + By =c, \nonumber
\end{align}
where $f(x)=\mathbb{E}_{\xi}[F(x,\xi)]$ is a nonconvex and smooth function; $\xi$ is a random vector;
$g(y)$ is nonsmooth and possibly nonconvex;
$x\in R^{p}$, $y\in R^{q}$, $A\in R^{d\times p}$, $B\in R^{d\times q}$ and $c \in R^d$.
The problem (1) is inspired by the structural risk minimization in machine learning \citep{vapnik2013nature}.
Here, the random vector $\xi$ obeys a fixed but unknown distribution, from which we are able to draw a set of i.i.d. samples.
In general, it is difficult to evaluate $\mathbb{E}_{\xi}[F(x,\xi)]$ exactly,
so we use the sample average approximation $\frac{1}{n}\sum_{i=1}^n F(x,\xi_i)$ to approximate it.
Throughout the paper, let $ f(x) \dot{=} \frac{1}{n}\sum_{i=1}^n F(x,\xi_i) = \frac{1}{n}\sum_{i=1}^n f_i(x)$
denote the average sum of many nonconvex and smooth component functions $f_i(x)$, $\forall i\in \{1,2,\cdots,n\}$.

Moreover, we propose three classes of nonconvex stochastic ADMM with variance reduction for the problem (1),
based on different reduced variance stochastic gradients.
Specifically, the first class called SVRG-ADMM uses a multi-stage scheme
to progressively reduce the variance of stochastic gradients.
The second called SAG-ADMM reduces the variance of stochastic gradient via additionally using the old gradients estimated in previous iteration.
The third called SAGA-ADMM is an extension of SAG-ADMM, which uses an unbiased stochastic gradient as the SVRG-ADMM.
In summary, our main contributions include three folds as follows:
\begin{itemize}
\item[1)] We propose three classes of the nonconvex stochastic ADMM with variance reduction, based on different reduced variance stochastic gradients.
\item[2)] We study the convergence of the proposed methods, and prove that these methods have
      the iteration complexity bound of $O(1/\epsilon)$ to obtain an $\epsilon$-stationary solution
      of the nonconvex problems. In particular, we provide a general framework
      to analyze the iteration complexity of the nonconvex stochastic ADMM with variance reduction.
\item[3)] Finally, some numerical experiments demonstrate the effectiveness of the proposed methods.
\end{itemize}

\subsection{Organization}
The paper is organized as follows:
In Section 2, we propose three classes of stochastic ADMM with variance reduction, based on different reduced variance
stochastic gradients.
In Section 3, we study the convergence and iterative complexity of the proposed methods.
Section 4 presents some numerical experiments, whose results back up the effectiveness of our methods.
In Section 5, we give some conclusions.
Most details of the theoretical analysis and proofs are relegated to the the following Appendices.

\subsection{Notations}
$\|\cdot\|$ denotes the Euclidean norm of a vector or the spectral norm of a matrix.
$\partial f$ is subgradient of the function $f$.
$Q\succ0$ implies the matrix $Q$ is positive definite. Let $\|x\|^2_Q = x^TQx$.
Let $A^+$ denote the generalized inverse of matrix $A$.
For a nonempty closed set $\mathcal {C}$, $\mbox{dist}(x,\mathcal {C})=\inf_{y\in \mathcal {C}}\|x-y\|$
denotes the distance from $x$ to $\mathcal {C}$.

\section{Stochastic ADMM methods for the Nonconvex Optimizations}

In this section, we study stochastic ADMM methods for solving the nonconvex problem (1).
First, we propose a simple nonconvex stochastic ADMM as a baseline, 
in which the variance of stochastic gradients is free.
However, it is difficult to guarantee the convergence of this simple stochastic ADMM, 
and only obtain a slow convergence rate.
Thus, we propose three classes of stochastic ADMM with variance reduction, 
based on different reduced variance stochastic gradients.

First, we review the standard ADMM for solving the problem (1), 
when $\xi$ is deterministic.
The augmented Lagrangian function of (1) is defined as
\begin{align}
 \mathcal {L}_{\rho}(x,y,\lambda) & = f(x) + g(y) - \langle\lambda, Ax+By-c\rangle + \frac{\rho}{2} \|Ax+By-c\|^2,
\end{align}
where $\lambda$ is a Lagrange multiplier, and $\rho$ is a penalty parameter.
At $t$-th iteration, the ADMM executes the following update:
\begin{align}
& y_{t+1} = \arg\min_y \mathcal {L}_{\rho}(x_t,y,\lambda_t) \\
& x_{t+1} = \arg\min_x \mathcal {L}_{\rho}(x,y_{t+1},\lambda_t) \\
& \lambda_{t+1} = \lambda_t - \rho(Ax_{t+1}+By_{t+1}-c).
\end{align}

When $\xi$ is a random variable, we can still update the variables $y$ and $\lambda$ by (3) and (5), respectively.
However, to update the variable $x$,
we need define an \emph{approximated} function of the form:
\begin{align}
 \hat{\mathcal {L}}_{\rho}(x,y,\lambda,v,\bar{x}) = f(\bar{x})+v^T(x-\bar{x}) + \frac{\eta}{2}\|x-\bar{x}\|_Q^2
 -\langle\lambda, Ax+By-c\rangle + \frac{\rho}{2} \|Ax+By-c\|^2,
\end{align}
where $\mathbb{E}[v] = \nabla f(\bar{x})$, $\eta>0$ and $Q \succ 0$.
By minimizing (6) on the variable $x$, we have
\begin{align}
 x \leftarrow (\eta Q + \rho A^TA)^{-1}\big(\eta Q \bar{x} - v - \rho A^T( By - c - \frac{\lambda}{\rho})\big). \nonumber
\end{align}
When $A^TA$ is large, computing inversion of $\eta Q + \rho A^TA$ is expensive.
To avoid it, we can use the inexact Uzawa method \citep{zhang2011unified} to
linearize the last term in (6), and choose $Q=(I-\frac{\rho}{\eta} A^TA)$.
When $Q=(I-\frac{\rho}{\eta} A^TA)$,
by minimizing (6) on the variable $x$, we have
\begin{align}
 x \leftarrow \bar{x} -\frac{1}{\eta}\big(v + \rho A^T(A\bar{x}+By-c-\frac{\lambda}{\rho})\big). \nonumber
\end{align}

\begin{algorithm}[htb]
   \caption{ S-ADMM for Nonconvex Optimization }
   \label{1}
\begin{algorithmic}[1]
   \STATE {\bfseries Input:} $T$, and $\rho>0$;
   \STATE {\bfseries Initialize:} $x_0$, $y_0$ and $\lambda_0$;
   \FOR {$t=0,1,\cdots,T-1$}
   \STATE{}  Uniformly randomly pick $i_t$ from $\{1,2,\cdots,n\}$;
   \STATE{}  $y_{t+1}=\arg\min_y \mathcal {L}_{\rho}(x_t,y,\lambda_t)$;
   \STATE{}  $x_{t+1}=\arg\min_x \hat{\mathcal {L}}_{\rho}(x,y_{t+1},\lambda_t,\nabla f_{i_t}(x_t),x_t)$;
   \STATE{}  $\lambda_{t+1} = \lambda_t-\rho(Ax_{t+1}+By_{t+1}-c)$;
   \ENDFOR
   \STATE {\bfseries Output:} Iterate $x$ and $y$ chosen uniformly random from $\{x_{t},y_{t}\}_{t=1}^{T}$.
\end{algorithmic}
\end{algorithm}

Like as the initial stochastic ADMM \citep{ouyang2013stochastic} for convex problems,
we propose a simple stochastic ADMM (S-ADMM) as a \emph{baseline} for the problem (1).
The algorithmic framework of the S-ADMM is given in Algorithm 2.
Though $\mathbb{E}[\nabla f_{i_t}(x)]= \nabla f(x)$, there exists the variance
$\mathbb{E}\|\nabla f_{i_t}(x)-\nabla f(x)\|^2$ in stochastic process.
To guarantee its convergence, we should choose a time-varying stepsize $1/\eta_{t}$ in (6), as in \citep{ouyang2013stochastic}.
However, as stochastic learning proceeds, the gradual decreasing of the stepsize $1/\eta_t$
generally leads to a slow convergence rate.
In the following, thus, we propose three classes of stochastic ADMM with
variance reduction for the problem (1), based on different reduced variance stochastic gradients.

\subsection{ Nonconvex SVRG-ADMM}

\begin{algorithm}[htb]
   \caption{ SVRG-ADMM for Nonconvex Optimization }
   \label{alg:2}
\begin{algorithmic}[1]
   \STATE {\bfseries Input:} epoch length $m$, $T$, $S=[T/m]$, $\rho>0$;
   \STATE {\bfseries Initialize:} $\tilde{x}^0=x_m^0$, $y_m^0$ and $\lambda_m^0$;
   \FOR {$s=0,1,\cdots,S-1$}
   \STATE{} $x_0^{s+1}=x_{m}^s$, $y_0^{s+1}=y_{m}^s$ and $\lambda_0^{s+1}=\lambda_{m}^s$;
   \STATE{} $\nabla f(\tilde{x}^s)=\frac{1}{n}\sum_{i=1}^n\nabla f_i(\tilde{x}^s)$;
   \FOR {$t=0,1,\cdots,m-1$}
   \STATE{} Uniformly randomly pick $i_t$ from $\{1,2,\cdots,n\}$;
   \STATE{} $y^{s+1}_{t+1}=\arg\min_y \mathcal {L}_{\rho}(x^{s+1}_t,y,\lambda_t^{s+1})$;
   \STATE{} $\hat{\nabla} f(x_{t}^{s+1}) = \nabla f_{i_t}(x_{t}^{s+1})-\nabla f_{i_t}(\tilde{x}^s)+\nabla f(\tilde{x}^s)$;
   \STATE{} $x^{s+1}_{t+1}=\arg\min_x \hat{\mathcal {L}}_{\rho}(x,y_{t+1}^{s+1},\lambda_{t}^{s+1},\hat{\nabla} f(x_{t}^{s+1}),x_{t}^{s+1})$;
   \STATE{} $\lambda_{t+1}^{s+1} = \lambda_{t}^{s+1}-\rho(Ax_{t+1}^{s+1}+By_{t+1}^{s+1}-c)$;
   \ENDFOR
   \STATE{} $\tilde{x}^{s+1}= x_m^{s+1}$;
   \ENDFOR
   \STATE {\bfseries Output:} Iterate $x$ and $y$ chosen uniformly random from $\{(x_{t}^s,y_{t}^s)_{t=1}^{m}\}_{s=1}^S$.
\end{algorithmic}
\end{algorithm}
In the subsection, we propose a nonconvex SVRG-ADMM, via using a multi-stage scheme
to progressively reduce the variance of stochastic gradients.
The framework of the SVRG-ADMM method is given in Algorithm 1.
Specifically, in Algorithm 2, the stochastic gradient $\hat{\nabla} f(x_{t}^{s+1})$ is unbiased, i.e.,
$\mathbb{E}[\hat{\nabla} f(x_{t}^{s+1})]=\nabla f(x_{t}^{s+1})$,
and its variance is progressively reduced by computing the gradients of all sample one time in each outer loop.
In the following, we give an upper bound of the variance of the stochastic gradient $\hat{\nabla} f(x_{t}^{s+1})$.
\begin{lemma}
 In Algorithm 2, set $\Delta^{s+1}_t=\hat{\nabla}f(x^{s+1}_t)-\nabla f(x^{s+1}_t)$,
 where $\hat{\nabla}f(x^{s+1}_t)=\nabla f_{i_t}(x^{s+1}_t)-\nabla f_{i_t}(\tilde{x}^{s})+\nabla f(\tilde{x}^{s})$,
 then the following inequality holds
 \begin{align}
  \mathbb{E}\|\Delta^{s+1}_t\|^2 \leq L^2\|x^{s+1}_t-\tilde{x}^{s}\|^2,
 \end{align}
 where $\mathbb{E}\|\Delta^{s+1}_t\|^2$ denotes the variance of stochastic gradient $\hat{\nabla}f(x^{s+1}_t)$.
\end{lemma}

A detailed proof of Lemma 1 is provided in \hyperref[app:lemma 1]{Appendix A}.
Lemma 1 shows that variance of the stochastic gradient $\hat{\nabla}f(x^{s+1}_t)$ has an upper bound $O(\|x_t^{s+1}-\tilde{x}^s\|^2)$.
Due to $\tilde{x}^s = x^s_m $, as number of iterations increases,
both $x_t^{s+1}$ and $\tilde{x}^s$ approach the same stationary point,
thus the variance of stochastic gradient vanishes.
Note that the stochastic ADMM for solving the nonconvex problem is difficult to converge to the global solution $x^*$,
so we bound the variance with $O(\|x_t^{s+1}-\tilde{x}^s\|^2)$ rather than the popular
$O(\|x-x^*\|^2)$ used in the convex problem.

\subsection{Nonconvex SAG-ADMM}

\begin{algorithm}[htb]
   \caption{ SAG-ADMM for Nonconvex Optimization }
   \label{3}
\begin{algorithmic}[1]
   \STATE {\bfseries Input:} $x_0 \in R^d$, $y_0\in R^q$, $z_i^0=x_0$ for $i\in \{1,2,\cdots,n\}$, number of iterations $T$;
   \STATE {\bfseries Initialize:} $\psi_0=\frac{1}{n}\sum_{i=1}^n\nabla f_i(z^0_i)$;
   \FOR {$t=0,1,\cdots,T-1$}
   \STATE{} Uniformly randomly pick $i_t,\ j_t$ from $\{1,2,\cdots,n\}$;
   \STATE{} $y_{t+1}=\arg\min_y \mathcal {L}_{\rho}(x_t,y,\lambda_t)$;
   \STATE{} $\hat{\nabla} f(x_{t}) = \frac{1}{n}\big(\nabla f_{i_t}(x_{t})-\nabla f_{i_t}(z^t_{i_t})\big)+\psi_t$;
   \STATE{} $x_{t+1}=\arg\min_x \hat{\mathcal {L}}_{\rho}(x,y_{t+1},\lambda_{t},\hat{\nabla} f(x_{t}),x_{t})$;
   \STATE{} $\lambda_{t+1} = \lambda_{t}-\rho(Ax_{t+1}+By_{t+1}-c)$;
   \STATE{} $z^{t+1}_{j_t}= x_t$ and $z_j^{t+1}=z^t_j$ for $j\neq j_t$;
   \STATE{} $\psi_{t+1}=\psi_t-\frac{1}{n}(\nabla f_{j_t}(z^t_{j_t})-\nabla f_{j_t}(z^{t+1}_{j_t}))$;
   \ENDFOR
   \STATE {\bfseries Output:} Iterate $x$ and $y$ chosen uniformly random from $\{x_{t},y_{t}\}_{t=1}^{T}$.
\end{algorithmic}
\end{algorithm}

In the subsection, we propose a nonconvex SAG-ADMM
by additionally using the old gradients estimated in the previous iteration.
The framework of the SAG-ADMM method is given in Algorithm 3.
In Algorithm 3, though the stochastic gradient $\hat{\nabla} f(x_{t})$ is biased, i.e.,
\begin{align}
 \mathbb{E}[\hat{\nabla} f(x_{t})]= \frac{1}{n}\nabla f(x_{t}) + (1-\frac{1}{n})\psi_t \neq \nabla f(x_{t}), \nonumber
\end{align}
its variance is progressively reduced.
In the following, we give an upper bound of variance of the stochastic gradient $\hat{\nabla} f(x_{t})$.

\begin{lemma}
 In Algorithm 3, set $\Delta_t=\hat{\nabla}f(x_t)-\nabla f(x_t)$,
 where $\hat{\nabla} f(x_{t}) = \frac{1}{n} \big(\nabla f_{i_t}(x_{t})-\nabla f_{i_t}(z^t_{i_t})\big)+\psi_t$, then the following inequality holds
 \begin{align}
  \mathbb{E}\|\Delta_t\|^2 \leq \big(1-\frac{1}{n}\big)^2\frac{L^2}{n}\sum_{i=1}^n \|x_t-z_i^t\|^2,
 \end{align}
 where $\psi_t=\frac{1}{n}\sum_{j=1}^n \nabla f_{j}(z^t_{j}) $, and $\mathbb{E}\|\Delta_t\|^2$
 denotes variance of the stochastic gradient $\hat{\nabla}f(x_t)$.
\end{lemma}

A detailed proof of Lemma 2 is provided in \hyperref[app:lemma 2]{Appendix B}.
Lemma 2 shows that the variance of the stochastic gradient $\hat{\nabla}f(x_t)$ has
an upper bound of $O(\big(1-\frac{1}{n}\big)^2\frac{1}{n}\sum_{i=1}^n \|x_t-z_i^t\|^2)$.
As the number of iteration increases,  both $x_t$ and the stored points $\{z^t\}_{i=1}^n$ approach the same stationary point,
so the variance of stochastic gradient progressively reduces.

\subsection{Nonconvex SAGA-ADMM}

\begin{algorithm}[htb]
   \caption{ SAGA-ADMM for Nonconvex Optimization }
   \label{4}
\begin{algorithmic}[1]
   \STATE {\bfseries Input:} $x_0 \in R^d$, $y_0\in R^q$, $z_i^0=x_0$ for $i\in \{1,2,\cdots,n\}$, number of iterations $T$;
   \STATE {\bfseries Initialize:} $\psi_0=\frac{1}{n}\sum_{i=1}^n\nabla f_i(z^0_i)$;
   \FOR {$t=0,1,\cdots,T-1$}
   \STATE{} Uniformly randomly pick $i_t,\ j_t$ from $\{1,2,\cdots,n\}$;
   \STATE{} $y_{t+1}=\arg\min_y \mathcal {L}_{\rho}(x_t,y,\lambda_t)$;
   \STATE{} $\hat{\nabla} f(x_{t}) = \nabla f_{i_t}(x_{t})-\nabla f_{i_t}(z^t_{i_t})+\psi_t$;
   \STATE{} $x_{t+1}=\arg\min_x \hat{\mathcal {L}}_{\rho}(x,y_{t+1},\lambda_{t},\hat{\nabla} f(x_{t}),x_{t})$;
   \STATE{} $\lambda_{t+1} = \lambda_{t}-\rho(Ax_{t+1}+By_{t+1}-c)$;
   \STATE{} $z^{t+1}_{j_t}= x_t$ and $z_j^{t+1}=z^t_j$ for $j\neq j_t$;
   \STATE{} $\psi_{t+1}=\psi_t-\frac{1}{n}(\nabla f_{j_t}(z^t_{j_t})-\nabla f_{j_t}(z^{t+1}_{j_t}))$;
   \ENDFOR
   \STATE {\bfseries Output:} Iterate $x$ and $y$ chosen uniformly random from $\{x_{t},y_{t}\}_{t=1}^{T}$.
\end{algorithmic}
\end{algorithm}

\begin{table}
  \centering
  \caption{Summary of the stochastic gradients used in the proposed methods. Note that $\nabla f(\tilde{x}^s)=\frac{1}{n}\sum_{i=1}^nf_{i}(\tilde{x}^s)$
  and $\psi_t=\frac{1}{n}\sum_{i=1}^nf_{i}(z^t_i)$. }\label{1}
  \begin{tabular}{|c|c|c|c|}
  \hline
  Methods & stochastic gradient & upper bound of variance & (un)biased  \\ \hline
   S-ADMM & $ \nabla f_{i_t}(x_{t})$  & unknown & unbiased  \\ \hline
   SVRG-ADMM & $ \nabla f_{i_t}(x_{t}^{s+1})-\nabla f_{i_t}(\tilde{x}^s)+\nabla f(\tilde{x}^s)$  & $L^2\|x^{s+1}_t-\tilde{x}^{s}\|^2$ & unbiased  \\ \hline
   SAG-ADMM & $\frac{1}{n}\big(\nabla f_{i_t}(x_{t})-\nabla f_{i_t}(z^t_{i_t})\big)+\psi_t$ & $\big(1-\frac{1}{n}\big)^2\frac{L^2}{n}\sum_{i=1}^n \|x_t-z_i^t\|^2$ &  biased  \\ \hline
   SAGA-ADMM &$\nabla f_{i_t}(x_{t})-\nabla f_{i_t}(z^t_{i_t})+\psi_t$ & $\frac{L^2}{n}\sum_{i=1}^n \|x_t-z_i^t\|^2$ & unbiased   \\
  \hline
  \end{tabular}
\label{tab:1}
\end{table}

In the subsection, we propose a nonconvex SAGA-ADMM,
which is an extension of the SAG-ADMM and uses an unbiased stochastic gradients as the SVRG-ADMM.
The framework of the SAGA-ADMM method is given in Algorithm 4.
In Algorithm 4, the stochastic gradient $\hat{\nabla} f(x_{t})$ is unbiased,
i.e., $\mathbb{E}[\hat{\nabla} f(x_{t})]=\nabla f(x_{t})$,
and its variance is progressively reduced via additionally also using the old gradients in the previous iterations.
Similarly, we give an upper bound of the variance of the stochastic gradient $\hat{\nabla} f(x_{t})$.

\begin{lemma}
 In Algorithm 4, set $\Delta_t=\hat{\nabla}f(x_t)-\nabla f(x_t)$,
 where $\hat{\nabla} f(x_{t}) = \nabla f_{i_t}(x_{t})-\nabla f_{i_t}(z^t_{i_t})+\psi_t$, then the following inequality holds
 \begin{align}
  \mathbb{E}\|\Delta_t\|^2 \leq \frac{L^2}{n}\sum_{i=1}^n \|x_t-z_i^t\|^2,
 \end{align}
 where $\psi_t=\frac{1}{n}\sum_{j=1}^n \nabla f_{j}(z^t_{j}) $, and $\mathbb{E}\|\Delta_t\|^2$
 denotes variance of the stochastic gradient $\hat{\nabla}f(x_t)$.
\end{lemma}

A detailed proof of Lemma 3 is provided in \hyperref[app:lemma 3]{Appendix C}.
Lemma 3 shows that the variance of the stochastic gradient $\hat{\nabla}f(x_t)$ has an upper bound $O(\frac{1}{n}\sum_{i=1}^n \|x_t-z_i^t\|^2)$.
Similarly, both $x_t$ and the stored points $\{z^t\}_{i=1}^n$ approach the same stationary point, as the number of iteration increases.
Thus, the variance of stochastic gradient progressively reduces. Note that the upper bound (9) loses a
coefficient $\big(1-\frac{1}{n}\big)^2$ to the upper bound (8), due to using a unbiased stochastic gradient in the SAGA-ADMM.

To further clarify the different of the proposed methods,  we summarize the stochastic gradients used in the proposed methods in Table \ref{tab:1}.
From Table \ref{tab:1}, we can find that the SAG-ADMM only uses the biased stochastic gradient,
while others use the unbiased stochastic gradient.
In particular, the SAG-ADMM can reduce faster the variance of stochastic gradient
than the SAGA-ADMM, at the expense of using a biased stochastic gradient.

\section{Convergence Analysis}

In the section, we analyze the convergence and iteration complexity of the proposed methods.
First, we give some mild assumptions regarding problem (1) as follows:
\begin{assumption}
For $\forall i \in \{1,2,\cdots,n\}$, the gradient of function $f_i$ is
Lipschitz continuous with the constant $L_i>0$, such that
\begin{align}
\|\nabla f_i(x_1)-\nabla f_i(x_2)\| \leq L_i \|x_1 - x_2\| \leq L \|x_1 - x_2\|, \ \forall x_1,x_2 \in R^p,
\end{align}
where $L=\max_i L_i$,
and this is equivalent to
\begin{align}
f_i(x_1) \leq f_i(x_2) + \nabla f_i(x_2)^T(x_1-x_2) + \frac{L}{2}\|x_1-x_2\|^2.
\end{align}
\end{assumption}

\begin{assumption}
$f(x) $ and $g(y)$ are all lower bounded, and denoting $f^*=\min_x f(x)$ and $g^*=\min_y g(y)$.
\end{assumption}

\begin{assumption}
$g(y)$ is a proper lower semi-continuous function.
\end{assumption}

\begin{assumption}
Matrix $A$ has full row rank.
\end{assumption}

In the Assumption 1, since $f(x)=\frac{1}{n}\sum_{i=1}^n f_i(x)$,
we have $\|\nabla f(x_1)-\nabla f(x_2)\| \leq L \|x_1 - x_2\|$ and
$f(x_1) \leq f(x_2) + \nabla f(x_2)^T(x_1-x_2) + \frac{L}{2}\|x_1-x_2\|^2, \ \forall x_1,x_2 \in R^p$.
Assumption 1 has been widely used in the convergence analysis
of nonconvex algorithms \citep{allen2016variance,reddi2016stochastic}.
Assumptions 2-3 have been used in study of ADMM for nonconvex problems \citep{jiang2016structured}.
Assumption 4 has been used in the convergence analysis of ADMM \citep{deng2016global}.

Throughout the paper, let $\sigma_{A}$ denote the smallest eigenvalues of matrix $AA^T$,
and let $\phi_{\min}$ and $\phi_{\max}$ denote the smallest and largest eigenvalues of positive matrix $Q$, respectively.
In the following, we define the $\epsilon$-stationary point of the nonconvex problem (1):

\begin{definition}
For $\epsilon>0$, the point $(x^*,y^*,\lambda^*)$ is said to be an $\epsilon$-stationary point of the problem (1) if it holds that
\begin{align}
\mathbb{E}\|\nabla f(x^*)-A^T\lambda^*\|^2 &\leq \epsilon, \\
\mathbb{E}\big[\mbox{dist}(B^T\lambda^*,\partial g(y^*))\big]^2 &\leq \epsilon,\\
\mathbb{E}\|Ax^*+By^*-c\|^2 &\leq \epsilon,
\end{align}
where $dist(y_0,\partial g(y)):=\inf \{\|y_0-z\|: \ z\in \partial g(y) \}$.
If $\epsilon=0$, the point $(x^*,y^*,\lambda^*)$ is said to be a stationary point of the problem (1).
\end{definition}
Note that combining the above inequalities (13-15)
is equivalent to $\mathbb{E}\big[ \mbox{dist}(0,\partial L(x^*,y^*,\lambda^*)) \big]^2 \leq \epsilon$, where
\begin{align}
   \partial L(x,y,\lambda) = \left [ \begin{matrix}
     \partial_x L(x,y,\lambda) \\
     \partial_y L(x,y,\lambda) \\
     \partial_{\lambda} L(x,y,\lambda)
 \end{matrix}
 \right ]. \nonumber
\end{align}
Next, based the above assumptions and definition, we study the convergence and iteration complexity of the proposed methods.
In particular, we provide a general framework to analyze the convergence
and iteration complexity of stochastic ADMM methods with variance reduction.
Specifically, the basic procedure is given as follows:
\begin{itemize}
  \item First, we design \textbf{a new sequence} based on the sequence generated from the algorithm. For example,
    we design the sequence $\{(\Psi^{s}_{t})_{t=1}^m\}_{s=1}^S$ in (15) for the SVRG-ADMM;
    the sequence $\{\Phi_{t}\}_{t=1}^T$ in (21) for the SAG-ADMM; and the sequence $\{\hat{\Phi}_{t}\}_{t=1}^T$ in (27) for the SAGA-ADMM.
  \item Second, we prove that the designed sequence is monotonically decreasing, and has a lower bound.
  \item Third, we define \textbf{a new variable} for the algorithm. For example, we define the variable $\theta^s_t$ in (19) for the SVRG-ADMM;
   the variable $\theta_t$ in (25) for the SAG-ADMM; and the variable $\hat{\theta}_t$ in (31) for the SAGA-ADMM. Then, we prove the variable has an upper bound,
   based on the above results.
  \item Finally, we prove that $\mathbb{E}\big[\mbox{dist}(0,\partial L(x,y,\lambda))\big]^2 $ is bounded by the above defined variable.
\end{itemize}

\subsection{Convergence Analysis of Nonconvex SVRG-ADMM}

In the subsection, we study the convergence and iteration complexity of the SVRG-ADMM.
First, given the sequence $\{(x^{s}_t,y^{s}_t,\lambda^{s}_t)_{t=1}^m\}_{s=1}^S$ generated by Algorithm 2,
then we define an useful sequence $\{(\Psi^{s}_{t})_{t=1}^m\}_{s=1}^S$ as follows:
\begin{align}
 \Psi^{s}_{t}
  =  \mathbb{E}\big[ \mathcal {L}_{\rho}(x^{s}_{t},y^{s}_{t},\lambda^{s}_{t})+ h^{s}_{t}(\|x^{s}_{t}-\tilde{x}^{s-1}\|^2+ \|x^{s}_{t-1}-\tilde{x}^{s-1}\|^2) + \frac{5(L^2+\eta^2\phi^2_{\max})}{\sigma_{A} \rho}\|x^{s}_{t}-x^{s}_{t-1}\|^2
   \big],
\end{align}
where the positive sequence $\{(h^s_t)_{t=1}^m \}_{s=1}^S$ satisfies the following formation (16).

Next, we consider three important lemmas: the first gives the upper bound of $\mathbb{E} \|\lambda^{s+1}_{t+1}-\lambda^{s+1}_{t}\|^2$;
the second demonstrates that the sequence $\{(\Psi^{s}_{t})_{t=1}^m\}_{s=1}^S$ is monotonically decreasing;
then the third give the lower bound of the sequence $\{(\Psi^{s}_{t})_{t=1}^m\}_{s=1}^S$.

\begin{lemma}
 Suppose the sequence $\{(x^{s}_t,y^{s}_t,\lambda^{s}_t)_{t=1}^m\}_{s=1}^S$ is generated by the Algorithm 2.
 The following inequality holds
 \begin{align}
   \mathbb{E}\|\lambda^{s+1}_{t+1}-\lambda^{s+1}_{t}\|^2 \leq & \frac{5L^2}{\sigma_{A}} \mathbb{E} \|x^{s+1}_{t}-\tilde{x}^{s}\|^2
   + \frac{5L^2}{\sigma_{A}}\|x^{s+1}_{t-1}-\tilde{x}^{s}\|^2
   + \frac{5\eta^2\phi^2_{\max}}{\sigma_{A}} \mathbb{E}\|x^{s+1}_{t+1}-x^{s+1}_t\|^2 \nonumber \\
  & + \frac{5(L^2+\eta^2\phi^2_{\max})}{\sigma_{A}}\|x^{s+1}_{t}-x^{s+1}_{t-1}\|^2, \nonumber
 \end{align}
 where $\sigma_{A}$ denotes the smallest eigenvalues of matrix $AA^T$,
  and $\phi_{\max}$ denotes the largest eigenvalues of positive matrix $Q$.
\end{lemma}
A detailed proof of Lemma 5 is provided in \hyperref[app:lemma 5]{Appendix D}.
Lemma 5 shows the upper bound of $\mathbb{E} \|\lambda^{s+1}_{t+1}-\lambda^{s+1}_{t}\|^2$.

\begin{lemma}
 Suppose that the sequence $\{(x^{s}_t,y^{s}_t,\lambda^{s}_t)_{t=1}^m\}_{s=1}^S$ is generated by Algorithm 2.
 Further suppose the positive sequence $\{(h_t^s)_{t=1}^m\}_{s=1}^S$ satisfies
 \begin{equation}
  h^{s}_t= \left\{
  \begin{aligned}
  & (2+\beta)h^{s}_{t+1}+ \frac{5L^2}{\sigma_{A} \rho}, \ 1 \leq t \leq m-1; \\
  & \frac{10L^2}{\sigma_{A}\rho}, \quad t=m,
  \end{aligned}
  \right.\end{equation}
 for all $s \in [S] $.
 Denoting
 \begin{equation}
 \Gamma^{s}_t =\left\{
 \begin{aligned}
   & \eta\phi_{\min} + \frac{\sigma_A\rho}{2} - \frac{L}{2} -\frac{5(2\eta^2\phi^2_{\max} + L^2)}{\sigma_{A}\rho}-(1+\frac{1}{\beta})h^{s}_{t+1}, \ 1 \leq t \leq m-1;  \\
   & \eta\phi_{\min} + \frac{\sigma_A\rho}{2} -\frac{L}{2}-\frac{5(2\eta^2\phi^2_{\max} + L^2)}{\sigma_{A}\rho} - h_1^{s+1}, \quad t=m
 \end{aligned}
 \right.\end{equation}
  and letting $\eta>0$, $\beta>0$ and $\rho>0$ be chosen such that
  \begin{align}
   \Gamma^{s}_t >0 , \ \forall t \in [m], \ \forall s \in [S], \nonumber
  \end{align}
 then the sequence $\{(\Psi^{s}_{t})_{t=1}^m\}_{s=1}^S$ is monotonically decreasing.
\end{lemma}

A detailed proof of Lemma 6 is provided in \hyperref[app:lemma 6]{Appendix E}.
Lemma 6 shows that the sequence $\{(\Psi^{s}_{t})_{t=1}^m\}_{s=1}^S$ is monotonically decreasing.
Next, we further clarify choosing the above parameters.
We first define a function $H(\eta)=\Gamma^s_{t}=\eta\phi_{\min} + \frac{\sigma_A\rho}{2}- \frac{L}{2}
- \frac{5(2\eta^2\phi^2_{\max}+L^2)}{\sigma_{A} \rho} -(1+\frac{1}{\beta})h^{s}_{t+1}$.
For the function $H(\eta)$, when $\eta=\frac{\sigma_A\rho \phi_{\min}}{20 \phi^2_{\max}}$,
the function $H(\eta)$ can reach the
largest value
\begin{align}
 H_{\max} = \frac{\sigma_A\rho}{40 \chi^2} + \frac{\sigma_A\rho}{2} - \frac{L}{2} - \frac{5L^2}{\sigma_A \rho} -(1+\frac{1}{\beta})h^{s}_{t+1}, \nonumber
\end{align}
where $\chi=\frac{\phi_{\max}}{\phi_{\min}}$ denotes the conditional number of matrix $Q$. Since $H(\eta)=\Gamma^s_{t}>0$,
we have $H_{\max}>0$. Considering $\sigma_A \rho >1$, then the parameter $\rho$ should satisfy the following inequality
\begin{align}
 \rho \geq \frac{40\chi^2}{(1+20\chi^2)\sigma_A}\big[\frac{L}{2} + 5L^2 + (1+\frac{1}{\beta})h^{s}_{t+1}\big].
\end{align}

\begin{lemma}
Suppose the sequence $\{(x^{s}_t,y^{s}_t,\lambda^{s}_t)_{t=1}^m\}_{s=1}^S$ is generated by the Algorithm 2,
and given the same conditions as in Lemma 6, the sequence $\{(\Psi^{s}_{t})_{t=1}^m\}_{s=1}^S$ has a lower bound.
\end{lemma}

A detailed proof of Lemma 7 is provided in \hyperref[app:lemma 7]{Appendix F}.
Next, based on the above lemmas, we analyze the convergence and iteration complexity of the SVRG-ADMM in the following.
We first defines an useful variable $\theta^s_t$ as follows:
 \begin{align}
  \theta^{s}_t = & \|x^{s}_{t}-\tilde{x}^{s-1}\|^2 + \|x^{s}_{t-1}-\tilde{x}^{s-1}\|^2
  + \|x^{s}_{t+1}-x^{s}_t\|^2 + \|x^{s}_{t}-x^{s}_{t-1}\|^2 .
 \end{align}

\begin{theorem}
 Suppose the sequence $\{(x^{s}_t,y^{s}_t,\lambda^{s}_t)_{t=1}^m\}_{s=1}^S$ is generated by the Algorithm 2.
 Denote
 $\kappa_1=3(L^2+\eta^2\phi^2_{\max})$, $\kappa_2=\frac{5(L^2+\eta^2\phi^2_{\max})}{\sigma_{A}\rho^2}$, $\kappa_3=\rho^2\|B\|_2^2\|A\|_2^2$,
 and $\tau = \min\big\{\gamma, \omega \big\}>0$, where $\gamma=\min_{t,s} \Gamma^s_t$ and $\omega = \frac{5L^2}{\sigma_{A} \rho}$.
 Letting
 \begin{align}
 mS = T = \frac{\max(\kappa_1,\kappa_2,\kappa_3)}{\tau \epsilon}(\Psi^{1}_{1}- \Psi^*),
 \end{align}
 where $\Psi^*$ is a lower bound of the sequence $\{(\Psi^{s}_{t})_{t=1}^m\}_{s=1}^S$,
 and denoting
 \begin{align}
  (\hat{t},\hat{s}) = \mathop{\arg\min}_{1 \leq t\leq m,\ 1 \leq s\leq S}\theta^{s}_t, \nonumber
 \end{align}
 then $(x_{\hat{t}}^{\hat{s}}, y_{\hat{t}}^{\hat{s}})$ is an $\epsilon$-stationary point of the problem (1).
\end{theorem}

A detailed proof of Theorem 8 is provided in \hyperref[app:theorem 8]{Appendix G}.
Theorem 8 shows that the SVRG-ADMM is convergent and has the iteration complexity
of $O(1/\epsilon)$ to reach an $\epsilon$-stationary point,
i.e., obtain a convergence rate $O(\frac{1}{T})$.
From Theorem 8, we can find that the SVRG-ADMM ensures its convergence
by progressively reducing the variance of stochastic gradients.

\subsection{Convergence Analysis of Nonconvex SAG-ADMM}

In the subsection, we study the convergence and iteration complexity of the SAG-ADMM.
First, given the sequence $\{x_t,y_t,\lambda_t\}_{t=1}^T$ generated by the Algorithm 3,
then we define an useful sequence $\{\Phi_{t}\}_{t=1}^T$ as follows
\begin{align}
 \Phi_t =   \mathbb{E} \big[ \mathcal {L}_{\rho}(x_{t},y_{t},\lambda_{t})
 + \big(1-\frac{1}{n}\big)^2 \frac{\alpha_{t}}{n}\sum_{i=1}^n(\|x_{t}-z^{t}_i\|^2 +\|x_{t-1}-z^{t-1}_i\|^2)
   + \frac{5(L^2+\eta^2\phi^2_{\max})}{\sigma_{A} \rho}\|x_{t}-x_{t-1}\|^2 \big],
\end{align}
where the positive sequence$\{\alpha_t\}_{t=1}^T$ satisfies the equality (22).

Next, we consider three important lemmas: the first gives 
the upper bound of $\mathbb{E}\|\lambda_{t+1}-\lambda_{t}\|^2$;
the second demonstrates that the sequence $\{\Phi_{t}\}_{t=1}^T$ is monotonically decreasing;
the third gives the lower bound of the sequence $\{\Phi_{t}\}_{t=1}^T$.

\begin{lemma}
 Suppose the sequence $\{x_t,y_t,\lambda_t\}_{t=1}^T$ is generated by the Algorithm 3, then the following inequality holds
 \begin{align}
  \mathbb{E}\|\lambda_{t+1}-\lambda_{t}\|^2 \leq & \big(1-\frac{1}{n}\big)^2\frac{5L^2}{\sigma_{A} n} \sum_{i=1}^n \mathbb{E} \|x_{t}-z^{t}_i\|^2
   + \big(1-\frac{1}{n}\big)^2\frac{5L^2}{\sigma_{A} n} \sum_{i=1}^n \|x_{t-1}-z^{t-1}_i\|^2 \nonumber \\
  & + \frac{5(\eta^2\phi^2_{\max})}{\sigma_{A}}\mathbb{E}\|x_{t+1}-x_t\|^2 + \frac{5(L^2+\eta^2\phi^2_{\max})}{\sigma_{A}}\|x_{t}-x_{t-1}\|^2. \nonumber
 \end{align}
\end{lemma}
A detailed proof of Lemma 9 is provided in \hyperref[app:lemma 9]{Appendix H}.
Lemma 9 shows the upper bound of $\mathbb{E}\|\lambda_{t+1}-\lambda_{t}\|^2$.
Next, we will prove that the sequence $\{\Phi_{t}\}_{t=1}^T$ is monotonically decreasing.

\begin{lemma}
 Suppose that the sequence $\{x_t,y_t,\lambda_t\}_{t=1}^T$ is generated by the Algorithm 3,
 and the positive sequence $\{\alpha_{t}\}_{t=1}^T$ satisfies
 \begin{align}
  \alpha_t=(2+\beta-\frac{1+\beta}{n})\alpha_{t+1} + \frac{5L^2}{\sigma_{A} \rho} + \frac{\vartheta L^2}{2}, \quad t=1,2,\cdots,T.
 \end{align}
 Denoting
 \begin{align}
\Gamma_{t} = \eta\phi_{\min} + \frac{\sigma_A\rho}{2}- \frac{L}{2} - \frac{1}{2\vartheta}-\frac{5(2\eta^2\phi^2_{\max} + L^2)}{\sigma_{A} \rho}
   -\big( 1-\frac{1}{n} \big)^2(1+\frac{1}{\beta}-\frac{1}{n\beta})\alpha_{t+1},
 \end{align}
 and Letting $\eta>0$, $\beta>0$, $\vartheta>0$ and $\rho>0$ be chosen such that $\Gamma_{t}>0$ for $t\geq 1$,
 then the sequence $\{\Phi_{t}\}_{t=1}^T$ is monotonically decreasing.
\end{lemma}

A detailed proof of Lemma 10 is provided in \hyperref[app:lemma 10]{Appendix I}.
Similarly, we further clarify how to choose the above parameters.
We first define a function $H(\eta)=\Gamma_{t}=\eta\phi_{\min} + \frac{\sigma_A\rho}{2}- \frac{L}{2} - \frac{1}{2\vartheta}- \frac{5(2\eta^2\phi^2_{\max}+L^2)}{\sigma_{A} \rho} -\big( 1-\frac{1}{n} \big)^2(1+\frac{1}{\beta}-\frac{1}{n\beta})\alpha_{t+1}$.
For the function $H(\eta)$, when $\eta=\frac{\sigma_A\rho \phi_{\min}}{20 \phi^2_{\max}}$, the function $H(\eta)$ can reach the
largest value
\begin{align}
 H_{\max} = \frac{\sigma_A\rho}{40 \chi^2} + \frac{\sigma_A\rho}{2}- \frac{L}{2} - \frac{1}{2\vartheta}- \frac{5L^2}{\sigma_A \rho} -
 \big( 1-\frac{1}{n} \big)^2(1+\frac{1}{\beta}-\frac{1}{n\beta})\alpha_{t+1}. \nonumber
\end{align}
Since $H(\eta)=\Gamma_{t}>0$, we have $H_{\max}>0$. Considering $\sigma_A \rho >1$,
then the parameter $\rho$ should satisfy the following inequality
\begin{align}
\rho \geq \frac{40\chi^2}{(1+20\chi^2)\sigma_A}\big[\frac{L}{2} + \frac{1}{2\vartheta} + 5L^2 + \big( 1-\frac{1}{n} \big)^2(1+\frac{1}{\beta}-\frac{1}{n\beta})\alpha_{1}\big].
\end{align}

\begin{lemma}
Suppose the sequence $\{x_t,y_t,\lambda_t\}_{t=1}^T$ is generated by the Algorithm 3.
Under the same conditions as in Lemma 10,
the sequence $\{\Phi_{t}\}_{t=1}^T$ has a lower bound.
\end{lemma}

A detailed proof of Lemma 11 is provided in \hyperref[app:lemma 11]{Appendix J}.
In the following, we will study the convergence and iteration complexity of the SAG-ADMM based on the above lemmas.
First, we defines an useful variable $\theta_t$ as follows:
 \begin{align}
  \theta_{t} = & \|x_{t+1}-x_t\|^2+\|x_t-x_{t-1}\|^2 + \big(1-\frac{1}{n}\big)^2\frac{1}{n}\sum_{i=1}^n \big( \|x_t-z^t_i\|^2 + \|x_{t-1}-z^{t-1}_i\|^2 \big).
 \end{align}

\begin{theorem}
 Suppose the sequence $\{x_t,y_t,\lambda_t\}_{t=1}^T$ is generated by Algorithm 3. Let
 $ \kappa_1=3(L^2+\eta^2\phi^2_{\max}), \ \kappa_2 = \frac{5(L^2+\eta^2\phi_{\max}^2)}{\sigma_{A}\rho^2}
 \ \mbox{and} \ \kappa_3=\rho^2\|B\|_2^2\|A\|_2^2$,
 and $\tau = \min\big\{ \gamma, \omega \big\}>0$,
 where $\gamma=\min_t\Gamma_t$ and $\omega=\min_t \big( 1-\frac{1}{n} \big)^2\big[(2+\beta-\frac{1+\beta}{n})\alpha_{t+1} + \frac{\vartheta L^2}{2}\big]$.
 Letting
 \begin{align}
 T =\frac{\max(\kappa_1,\kappa_2,\kappa_3)}{\epsilon\tau}(\Phi_{1}- \Phi^*),
 \end{align}
 where $\Phi^*$ denotes a lower bound of the sequence $\{\Phi_{t}\}_{t=1}^T$,
 and denoting
 \begin{align}
   \hat{t} = \mathop{\arg\min}_{1\leq t \leq T }\theta_{t}, \nonumber
 \end{align}
 then $(x_{\hat{t}},y_{\hat{t}})$ is an $\epsilon$-stationary point of the problem (1).
\end{theorem}

A detailed proof of Theorem 12 is provided in \hyperref[app:theorem 12]{Appendix K}.
Theorem 12 shows that the SAG-ADMM method
has the iteration complexity of $O(1/\epsilon)$ to reach an $\epsilon$-stationary point,
i.e., obtain a convergence rate $O(\frac{1}{T})$.

\subsection{Convergence Analysis of Nonconvex SAGA-ADMM}
In the subsection, we study the convergence and iteration complexity of the SAGA-ADMM.
Similarly, we first define an useful sequence $\{\hat{\Phi}_{t}\}_{t=1}^T$ as follows:
\begin{align}
 \hat{\Phi}_{t} =  \mathbb{E} \big[ \mathcal {L}_{\rho}(x_{t},y_{t},\lambda_{t})
 +\frac{\alpha_{t}}{n}\sum_{i=1}^n(\|x_{t}-z^{t}_i\|^2 +\|x_{t-1}-z^{t-1}_i\|^2)
  + \frac{5(L^2+\eta^2\phi^2_{\max})}{\sigma_{A} \rho}\|x_{t}-x_{t-1}\|^2 \big],
\end{align}
where the positive sequence$\{\alpha_t\}_{t=1}^T$ satisfies the equality (28).

Similarly, we consider three important lemmas: the first gives the upper bound of $\mathbb{E}\|\lambda_{t+1}-\lambda_{t}\|^2$;
the second demonstrates that the sequence $\{\hat{\Phi}_{t}\}_{t=1}^T$ is monotonically decreasing;
the third shows the lower bound of the sequence $\{\hat{\Phi}_{t}\}_{t=1}^T$.

\begin{lemma}
 Suppose the sequence $\{x_t,y_t,\lambda_t\}_{t=1}^T$ is generated by the Algorithm 4. The following inequality holds
 \begin{align}
  \mathbb{E}\|\lambda_{t+1}-\lambda_{t}\|^2 \leq & \frac{5L^2}{\sigma_{A} n} \sum_{i=1}^n \mathbb{E} \|x_{t}-z^{t}_i\|^2
   + \frac{5L^2}{\sigma_{A} n} \sum_{i=1}^n \|x_{t-1}-z^{t-1}_i\|^2
  + \frac{5(\eta^2\phi^2_{\max})}{\sigma_{A}}\mathbb{E}\|x_{t+1}-x_t\|^2 \nonumber \\
  & + \frac{5(L^2+\eta^2\phi^2_{\max})}{\sigma_{A}}\|x_{t}-x_{t-1}\|^2. \nonumber
 \end{align}
\end{lemma}
A detailed proof of lemma 13 is provided in \hyperref[app:lemma 13]{Appendix L}.
Lemma 13 shows that $\mathbb{E}\|\lambda_{t+1}-\lambda_{t}\|^2$ has an upper bound.

\begin{lemma}
 Suppose that the sequence $\{x_t,y_t,\lambda_t\}_{t=1}^T$ is generated by the Algorithm 4,
 and the positive sequence $\{\alpha_{t}\}_{t=1}^T$ satisfy
 \begin{align}
  \alpha_t=(2+\beta-\frac{1+\beta}{n})\alpha_{t+1}+\frac{5L^2}{\sigma_{A} \rho}, \quad t=1,2,\cdots,T.
 \end{align}
 Denoting
 \begin{align}
\Gamma_{t} = \eta\phi_{\min} + \frac{\sigma_A\rho}{2} - \frac{L}{2} - \frac{5(2\eta^2\phi^2_{\max}+L^2)}{\sigma_{A} \rho} -(1+\frac{1}{\beta}-\frac{1}{n\beta})\alpha_{t+1},
 \end{align}
 and Letting $\eta>0$, $\beta>0$ and $\rho>0$ be chosen such that $\Gamma_{t}>0$ for $t\geq 1$,
 then the sequence $\{\hat{\Phi}_t\}$ is monotonically decreasing.
\end{lemma}

A detailed proof of Lemma 14 is provided in \hyperref[app:lemma 14]{Appendix M}.
Lemma 14 shows that the sequence $\{\hat{\Phi}_t\}$ is monotonically decreasing.
Next, we further clarify how to choose the above parameters used in Lemma 14.
First, we define a function $H(\eta)=\Gamma_{t}=\eta\phi_{\min} + \frac{\sigma_A\rho}{2}- \frac{L}{2} - \frac{5(2\eta^2\phi^2_{\max}+L^2)}{\sigma_{A} \rho} -(1+\frac{1}{\beta}-\frac{1}{n\beta})\alpha_{t+1}$.
For the function $H(\eta)$, when $\eta=\frac{\sigma_A\rho \phi_{\min}}{20 \phi^2_{\max}}$, the function $H(\eta)$ can reach the
largest value
\begin{align}
 H_{\max} = \frac{\sigma_A\rho}{40 \chi^2} + \frac{\sigma_A\rho}{2}- \frac{L}{2} - \frac{5L^2}{\sigma_A \rho} -(1+\frac{1}{\beta}-\frac{1}{n\beta})\alpha_{t+1}. \nonumber
\end{align}
Since $H(\eta)=\Gamma_{t}>0$,
we have $H_{\max}>0$. Considering $\sigma_A \rho >1$, then the parameter $\rho$ should satisfy the following inequality
\begin{align}
\rho \geq \frac{40\chi^2}{(1+20\chi^2)\sigma_A}\big[\frac{L}{2} + 5L^2 + (1+\frac{1}{\beta}-\frac{1}{n\beta})\alpha_{1}\big].
\end{align}

\begin{lemma}
Suppose the sequence $\{x_t,y_t,\lambda_t\}_{t=1}^T$ is generated by the Algorithm 4. Under the same conditions as in Lemma 14,
the sequence $\hat{\Phi}_{t}$ has a lower bound.
\end{lemma}

The proof of Lemma 15 can follow the proof of Lemma 11.
In the following, we will study the convergence and iteration complexity of the SAGA-ADMM based on the above lemmas.
First, we defines an useful variable $\hat{\theta}_t$ as follows:
 \begin{align}
  \hat{\theta}_{t} = & \|x_{t+1}-x_t\|^2+\|x_t-x_{t-1}\|^2 + \frac{1}{n}\sum_{i=1}^n \big( \|x_t-z^t_i\|^2 + \|x_{t-1}-z^{t-1}_i\|^2 \big).
 \end{align}

\begin{theorem}
 Suppose the sequence $\{x_t,y_t,\lambda_t\}_{t=1}^T$ is generated by the Algorithm 4. Let
 $ \kappa_1=3(L^2+\eta^2\phi^2_{\max}), \ \kappa_2 = \frac{5(L^2+\eta^2\phi_{\max}^2)}{\sigma_{A}\rho^2}
 \ \mbox{and} \ \kappa_3=\rho^2\|B\|_2^2\|A\|_2^2$,
 and $\tau = \min\big\{ \gamma, \omega \big\}>0$,
 where $\gamma=\min_t\Gamma_t$ and $\omega=\min_t (2+\beta-\frac{1+\beta}{n})\alpha_{t+1}$.
 Letting
 \begin{align}
 T =\frac{\max(\kappa_1,\kappa_2,\kappa_3)}{\epsilon\tau}(\hat{\Phi}_{1} - \hat{\Phi}^*),
 \end{align}
 where $\hat{\Phi}^*$ is a lower bound of the sequence $\{\hat{\Phi}_{t}\}_{t=1}^T$,
 and denoting
 \begin{align}
   \tilde{t} = \mathop{\arg\min}_{1\leq t \leq T }\hat{\theta}_{t}, \nonumber
 \end{align}
 then $(x_{\tilde{t}},y_{\tilde{t}})$ is an $\epsilon$-stationary point of the problem (1).
\end{theorem}

A detailed proof of Theorem 16 is provided in \hyperref[app:theorem 16]{Appendix N}.
Theorem 16 shows that the SAGA-ADMM method has the iteration complexity of $O(1/\epsilon)$ to reach an $\epsilon$-stationary point, i.e.,
obtain a convergence rate $O(\frac{1}{T})$.

\subsection{Convergence Analysis for Nonconvex S-ADMM}
In the subsection, we will prove that the S-ADMM,
in which the variance of stochastic gradients is free, is divergent under some conditions.

\begin{theorem}
In Algorithm 1, given a constant stepsize parameter $\eta$, and let $\delta>0$ denote a constant.
Suppose the variance of stochastic gradients satisfy $ \mathbb{E}\|\nabla f_{i_t}(x)-\nabla f(x)\|^2 \geq \delta^2$.
If $\delta \geq 2(L+\eta\phi_{\max})\epsilon$, the S-ADMM will be divergent.
\end{theorem}

A detailed proof of Theorem 17 is provided in \hyperref[app:theorem 17]{Appendix O}.
Theorem 17 shows that when given a constant parameter $\eta$, the S-ADMM may be divergent.
In other words, the variance of stochastic algorithms
easily leads to the iteration points jumping from the neighbourhood
of a stationary point to that of another stationary point of the nonconvex problems.
Thus, we should consider controlling variance of the stochastic gradients,
when design the stochastic ADMM for the non-convex optimizations.

\section{Experiments}

In this section, we will execute some numerical experiments to demonstrate
the performances of the proposed methods for the nonconvex optimizations.
In the following, all algorithms are implemented in MATLAB,
and experiments are performed on a PC with an Intel i7-4770 CPU and 16GB memory.

\subsection{Experimental Setups}
In the experiments, we focus on the binary classification with incorporating the correlations between features.
Given a set of straining samples $\{(a_i,b_i)\}_{i=1}^n$,
where $a_i\in R^d$, $b_i \in \{-1,+1\}$, $\forall i \in \{1,2,\cdots,n\}$,
then we solve the following noncovex robust graph-guided models, as the graph-guided fused lasso \citep{kim2009multivariate}, 
\begin{align}
 \min_x \frac{1}{n}\sum_{i=1}^n f_i(x) + \lambda_1\|Ax\|_1 + \frac{\lambda_2}{2}\|x\|_2^2,
\end{align}
where $f_i(x)=\frac{1}{1+\exp(b_i a_i^Tx)}$ denotes the sigmoid loss function,
which is nonconvex and smooth;
$\lambda_1$ and $\lambda_2$ are positive regularization parameters.
Specifically, let $A=[G;I]$, where the matrix $G$ is obtained 
by sparse inverse covariance matrix estimation \citep{friedman2008sparse,hsieh2014quic}.
In order to satisfy the ADMM formulation, we can introduce an additional variable $y$ and 
rewrite the problem (33) as follows:
\begin{align}
 & \min_{x,y} f(x) + g(y) \nonumber \\
 & \mbox{s.t.} \ Ax-y=0, \nonumber
\end{align}
where $f(x)=\frac{1}{n}\sum_{i=1}^n f_i(x)+ \frac{\lambda_2}{2}\|x\|_2^2$ and $g(y)=\lambda_1\|y\|_1$.

\begin{table}
  \centering
  \caption{Summary of data sets and regularization parameters used in our experiments.}\label{1}
  \begin{tabular}{c|c|c|c|c}
  \hline
  data sets & number of samples & dimensionality & $\lambda_1$ & $\lambda_2$ \\ \hline
  \emph{a9a}   & 32,561 & 123 & $ 10^{-4}$ & $ 1.2 \times 10^{-4}$\\
  \emph{covertype} & 581,012 & 54 & $ 10^{-4}$ & $ 10^{-6}$ \\
  \emph{mnist8m} & 1,100,000 & 784 & $ 10^{-3}$ & $ 1.2 \times 10^{-3}$\\
  \hline
  \end{tabular}
 \label{tab:2}
\end{table}

In the experiments, we use three publicly available datasets\footnote{\emph{a9a}, \emph{covertype}, and \emph{mnist8m}
are from the LIBSVM website (www.csie.ntu.edu.tw/~cjlin/libsvmtools/datasets/).}, which are summarized in Table \ref{tab:2}.
For each dataset, we use half of the samples as training data, while use the rest as testing data.
Note that we only consider classifying the first class in the dataset \emph{mnist8m}.
In the algorithms, we choose the initial solution $x_0=zeros(d,1)$ and $\lambda_0=A^+ \nabla f(x_0)$.
At the same time, we fix the parameters $\eta=2$ and $\rho=6$, and set $Q=I$.
In particular, we consider two cases in Algorithm 1: the \emph{S-ADMM} with a time-varying stepsize parameter
$\eta=2\sqrt{t}$; the \emph{S-ADMM-F} with a fixed parameter $\eta=2$.
In Table \ref{tab:2}, we also provide some regularization parameters used in our experiments.
In the SVRG-ADMM algorithm, we choose $m=n$.
Finally, all experimental results are averaged over 10 repetitions.

\begin{figure*}[htbp]
\centering
\subfigure[\emph{a9a}]{\includegraphics[width=0.3\textwidth]{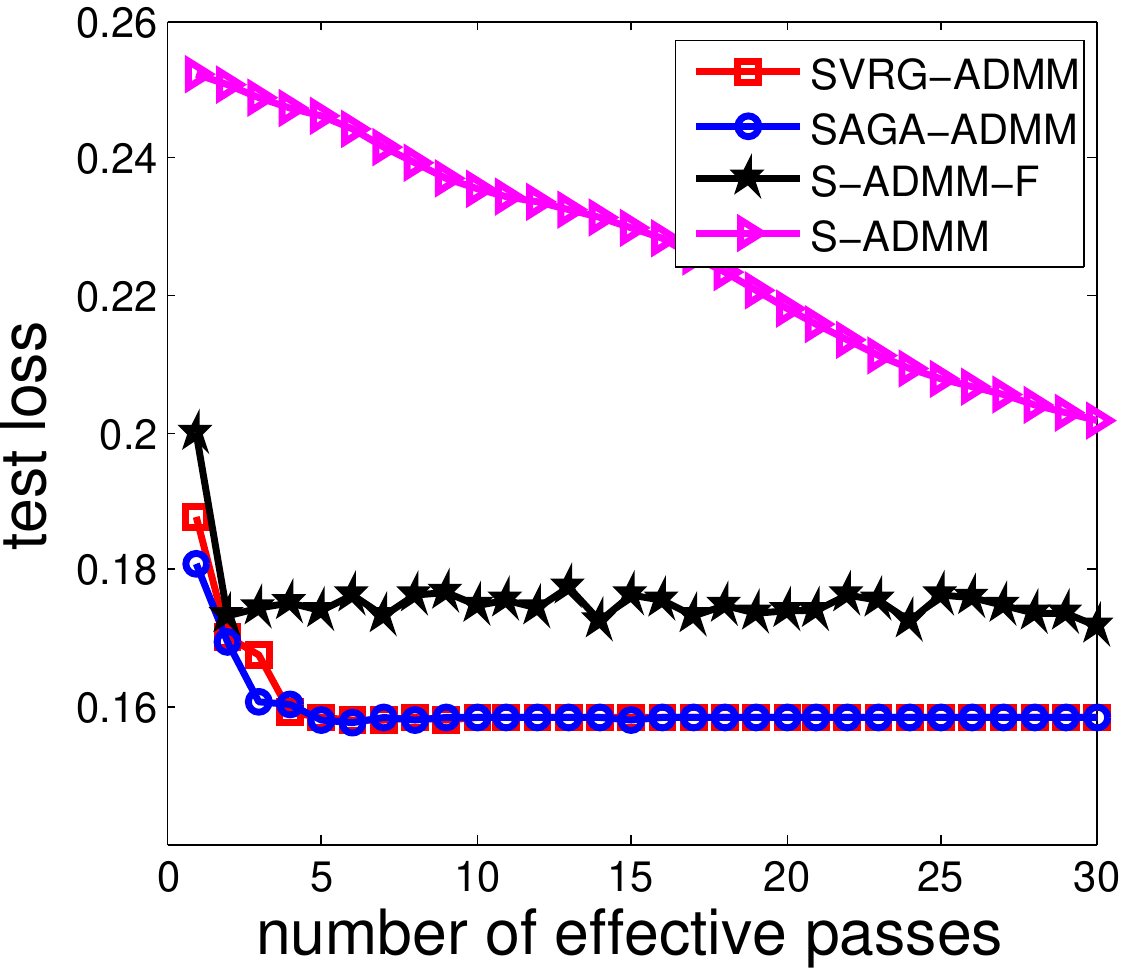}}
\subfigure[\emph{covertype}]{\includegraphics[width=0.3\textwidth]{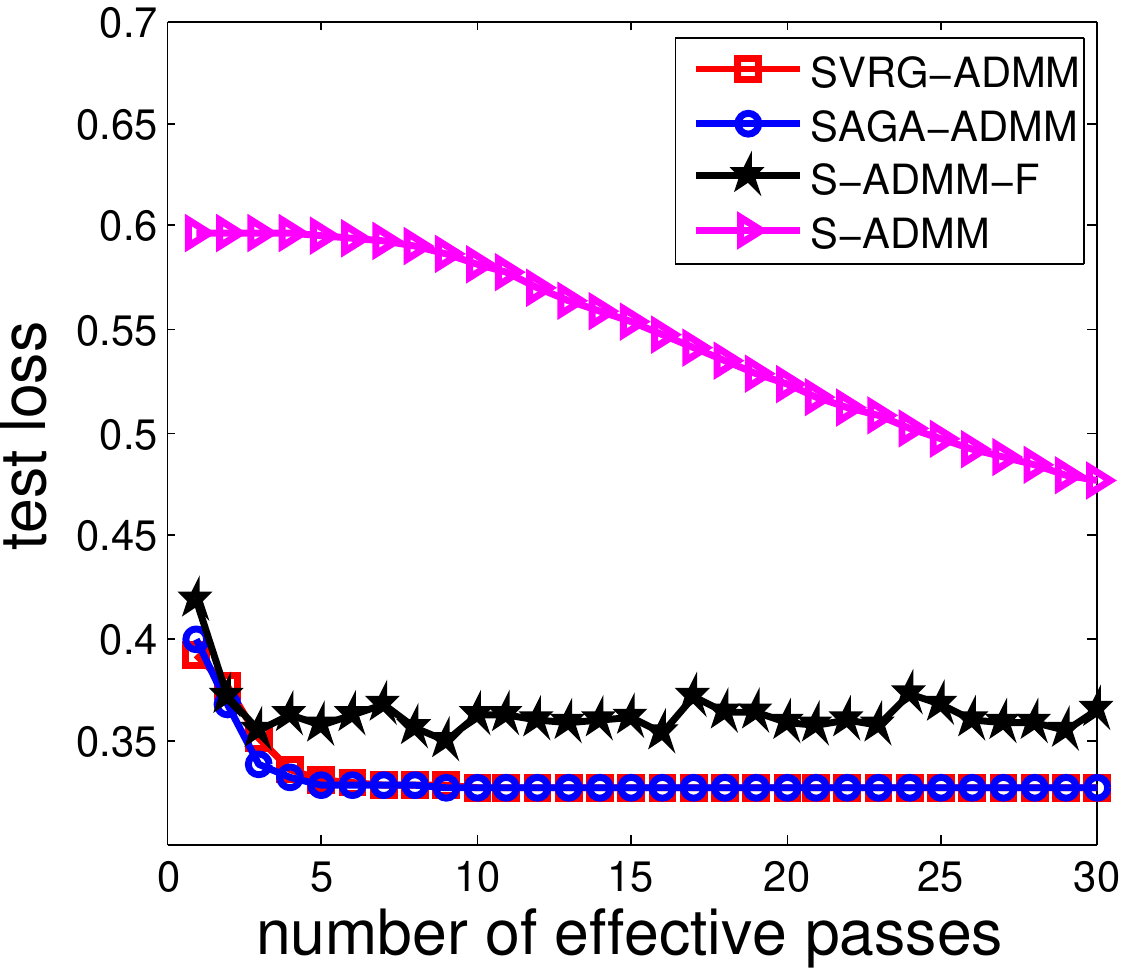}}
\subfigure[\emph{mnist8m}]{\includegraphics[width=0.3\textwidth]{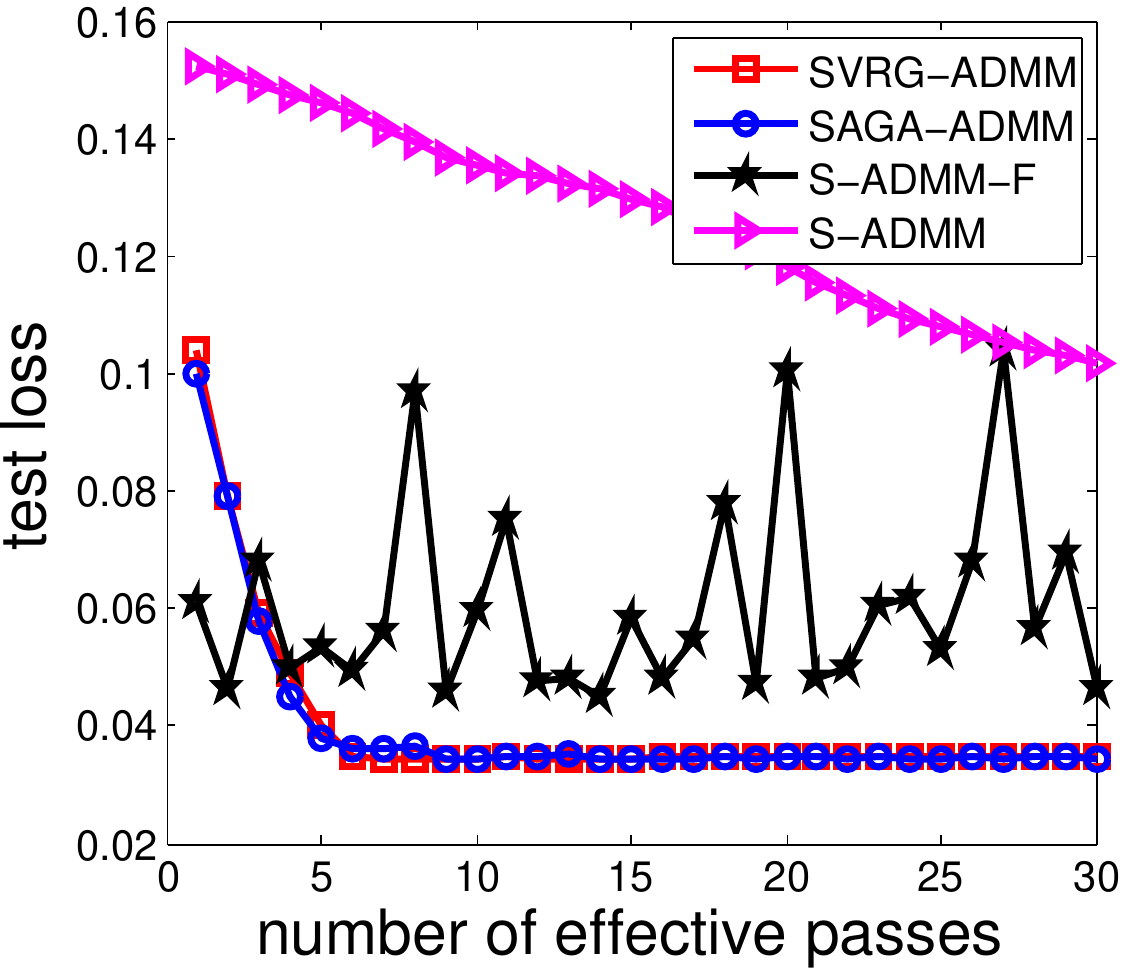}}
\caption{Test loss \emph{versus} number of effective passes on the \emph{nonconvex} graph-guided model.} \label{1}
\vspace{-1em}
\end{figure*}

\begin{figure*}[htbp]
\centering
\subfigure[\emph{a9a}]{\includegraphics[width=0.3\textwidth]{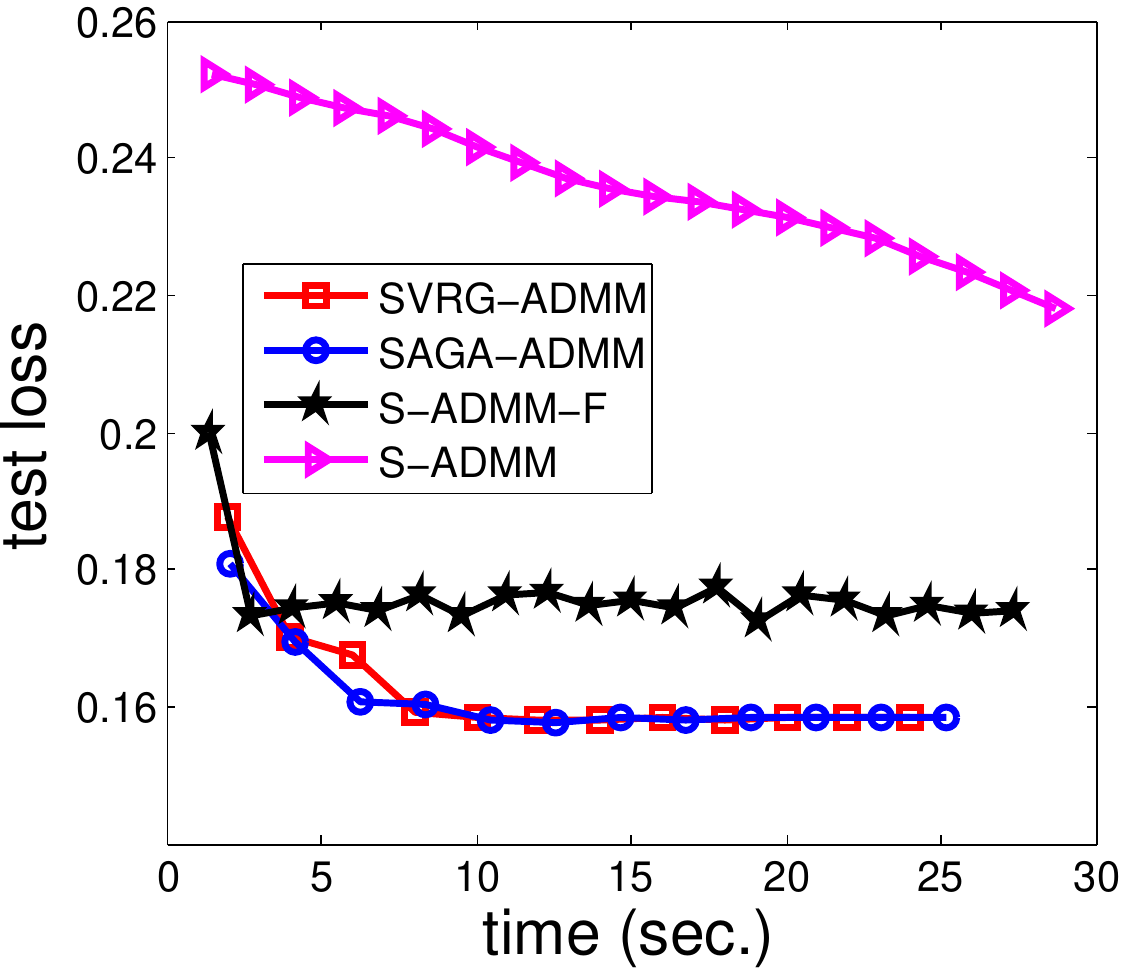}}
\subfigure[\emph{covertype} ]{\includegraphics[width=0.3\textwidth]{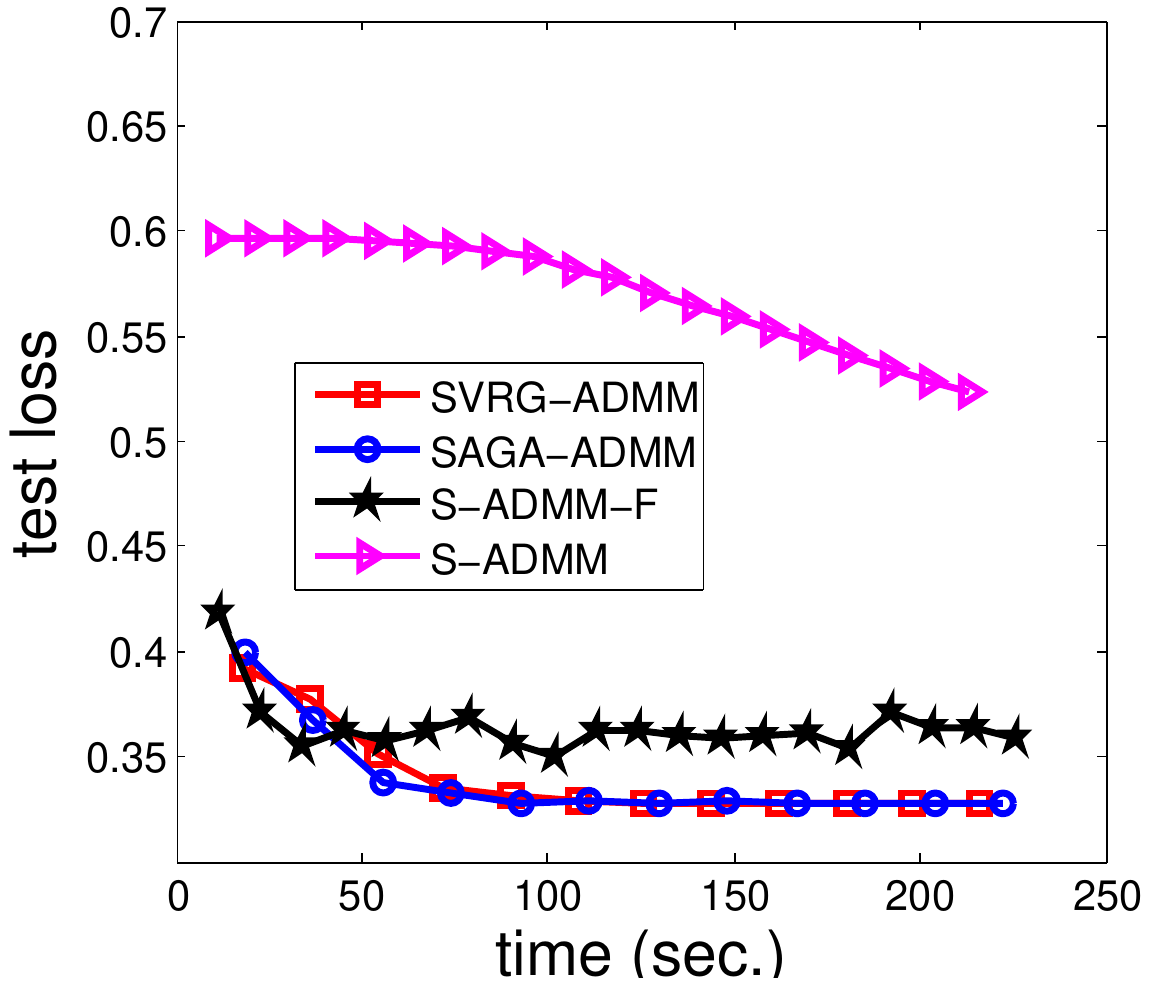}}
\subfigure[\emph{mnist8m}]{\includegraphics[width=0.3\textwidth]{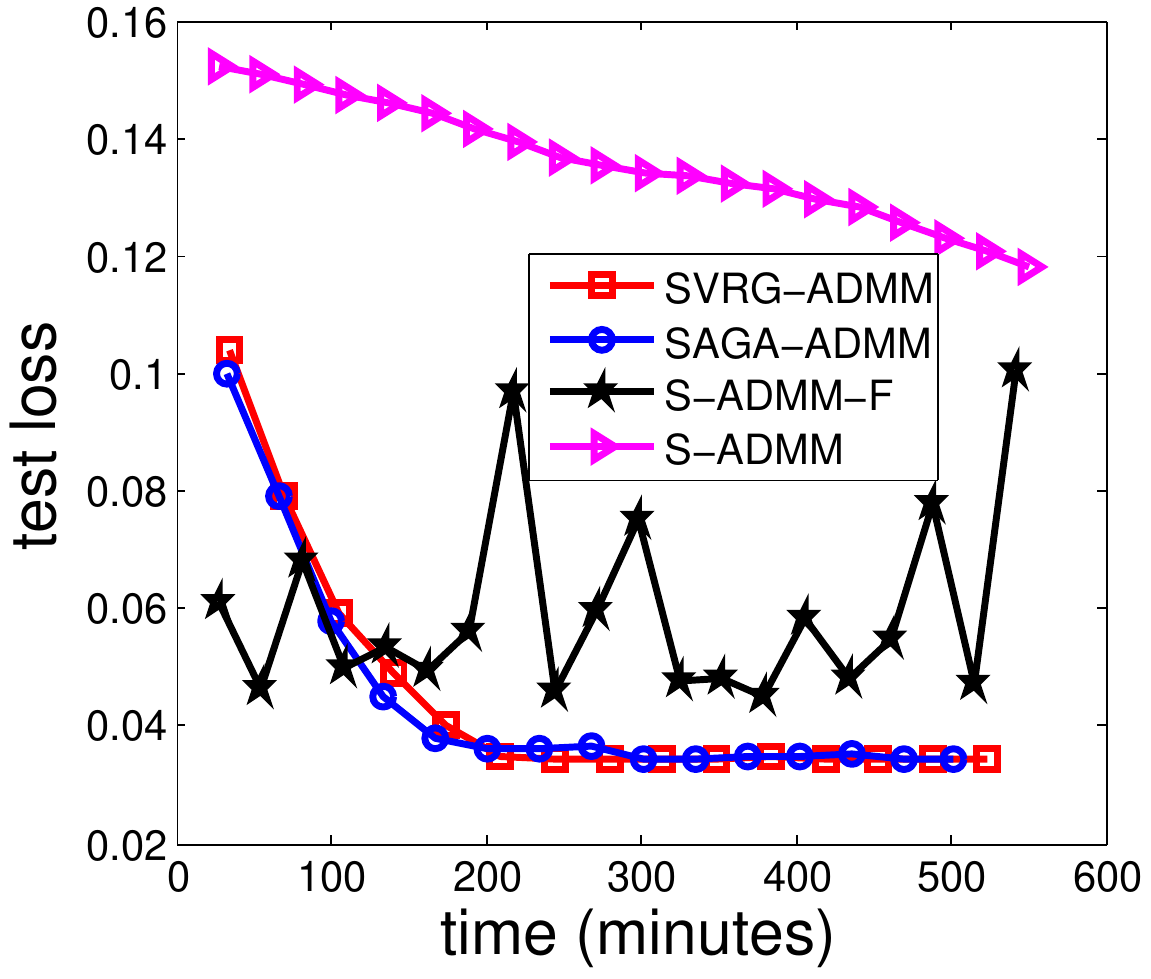}}
\caption{Test loss \emph{versus} time on the \emph{nonconvex} graph-guided model.} \label{1}
\vspace{-1em}
\end{figure*}

\begin{figure*}[htbp]
\centering
\subfigure[\emph{a9a}]{\includegraphics[width=0.3\textwidth]{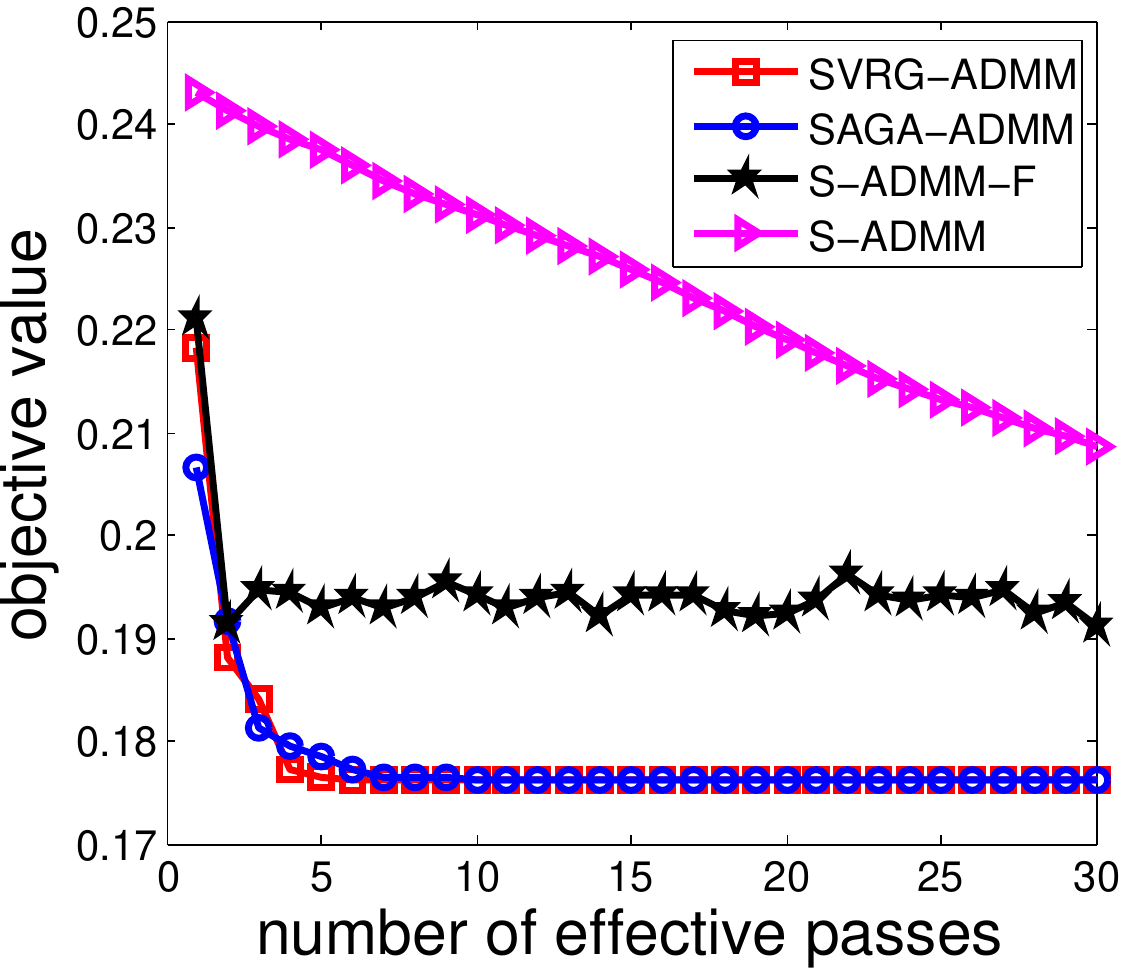}}
\subfigure[\emph{covertype}]{\includegraphics[width=0.3\textwidth]{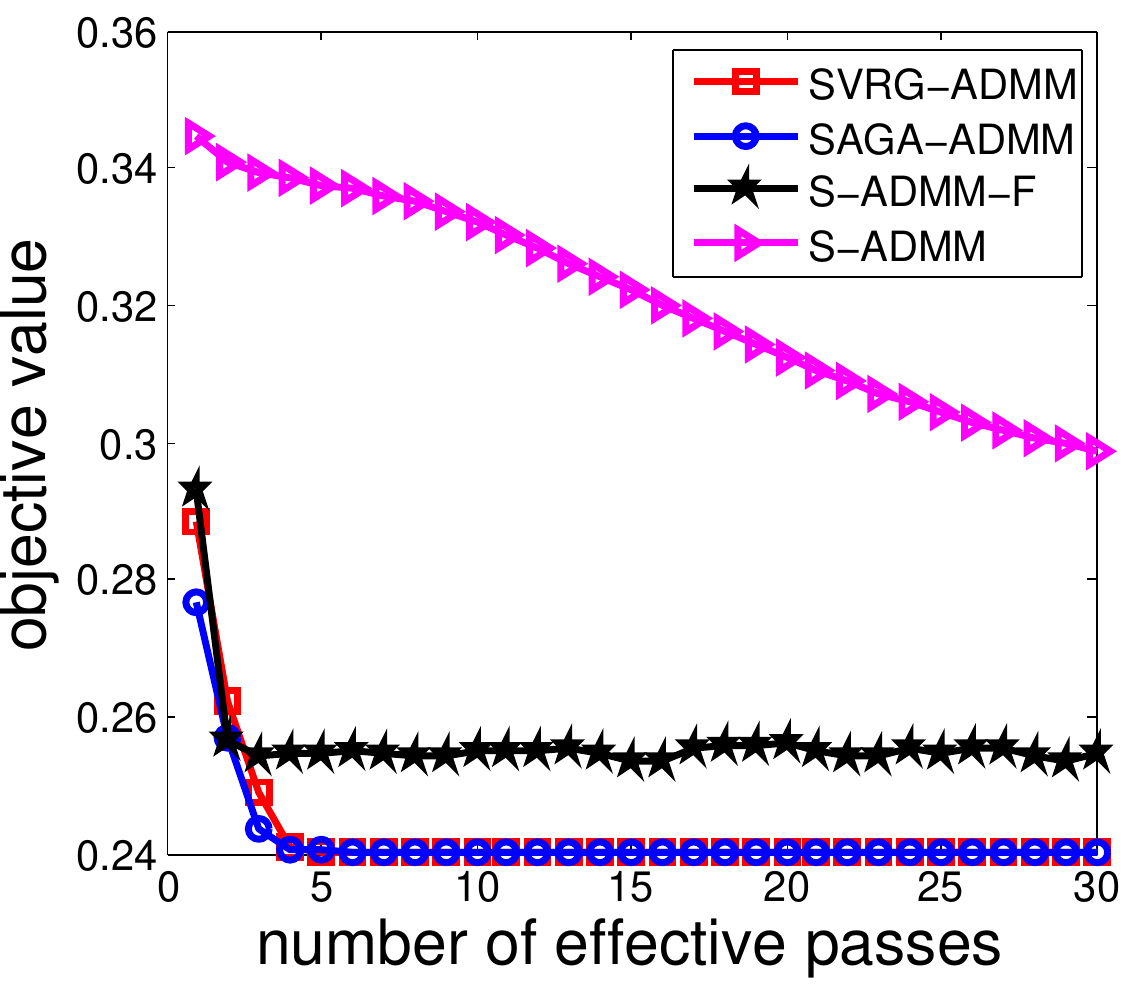}}
\subfigure[\emph{mnist8m}]{\includegraphics[width=0.3\textwidth]{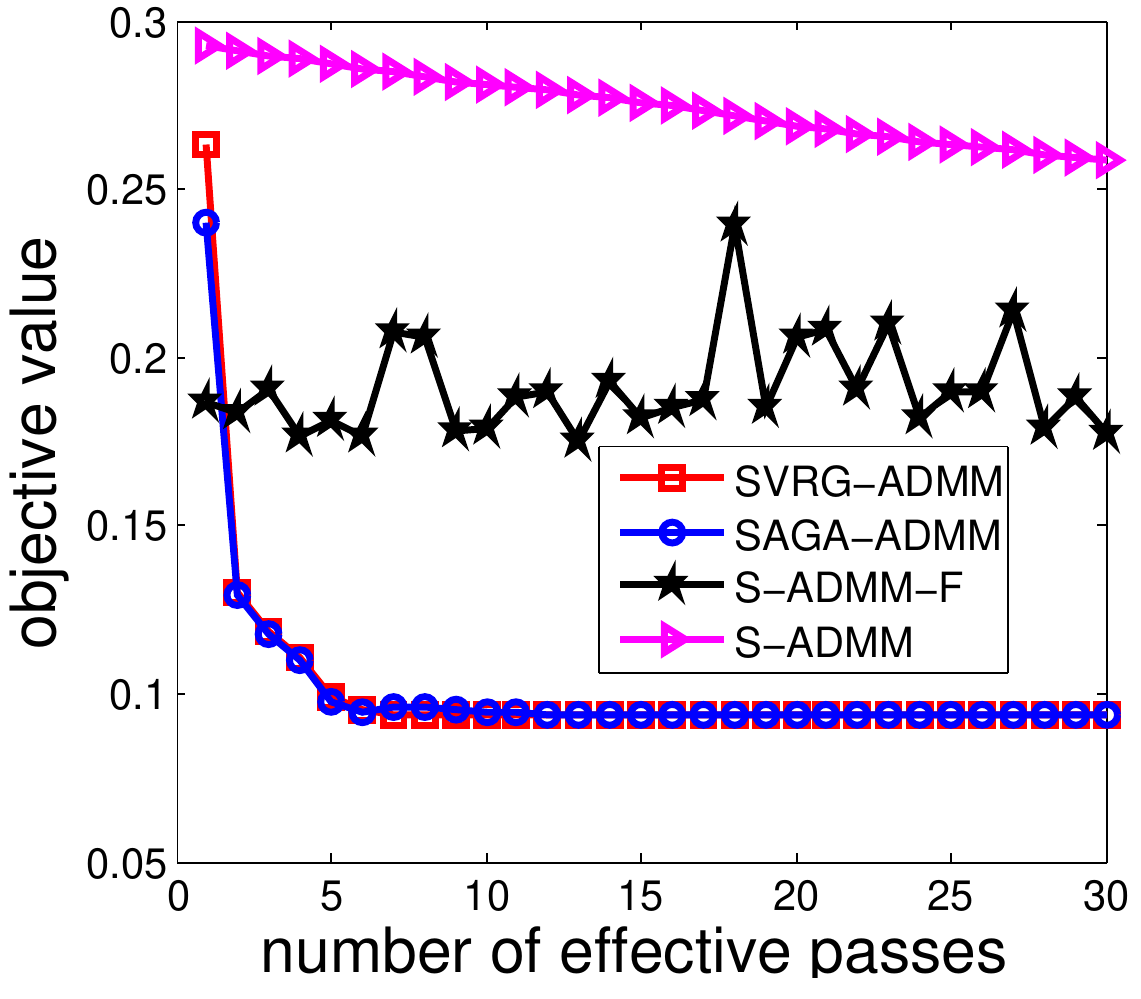}}
\caption{Objective value \emph{versus} number of effective passes on the \emph{nonconvex} graph-guided model.} \label{1}
\vspace{-1em}
\end{figure*}

\begin{figure*}[htbp]
\centering
\subfigure[\emph{a9a}]{\includegraphics[width=0.3\textwidth]{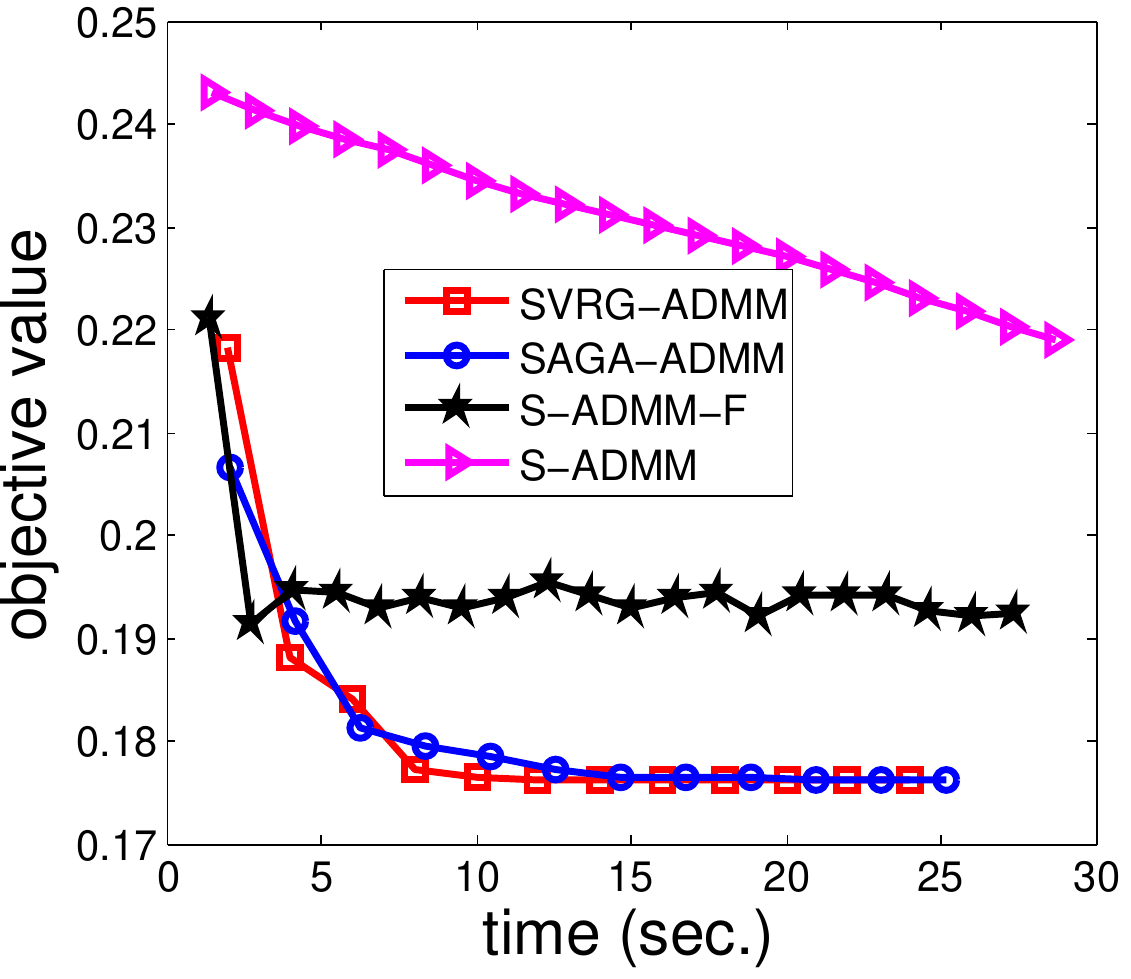}}
\subfigure[\emph{covertype}]{\includegraphics[width=0.3\textwidth]{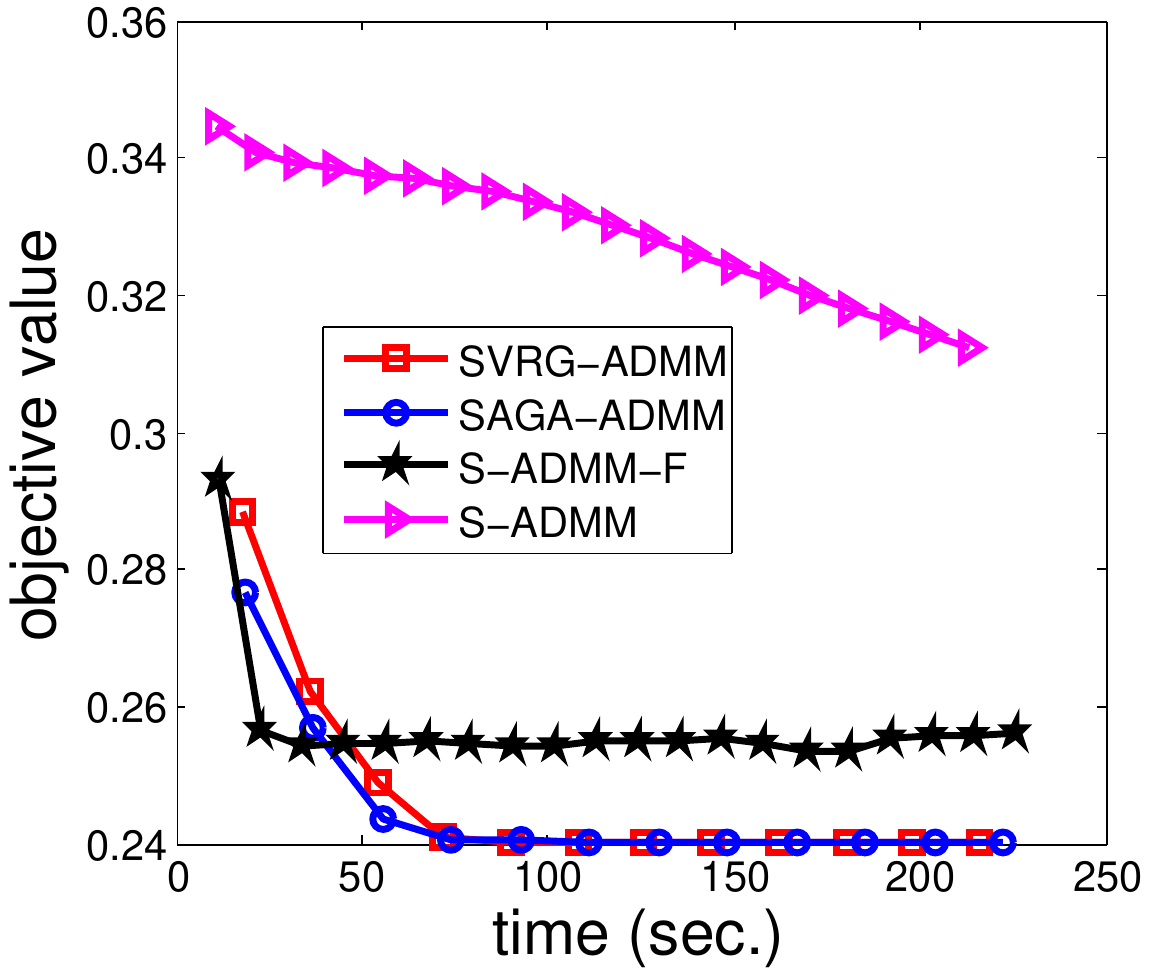}}
\subfigure[\emph{mnist8m}]{\includegraphics[width=0.3\textwidth]{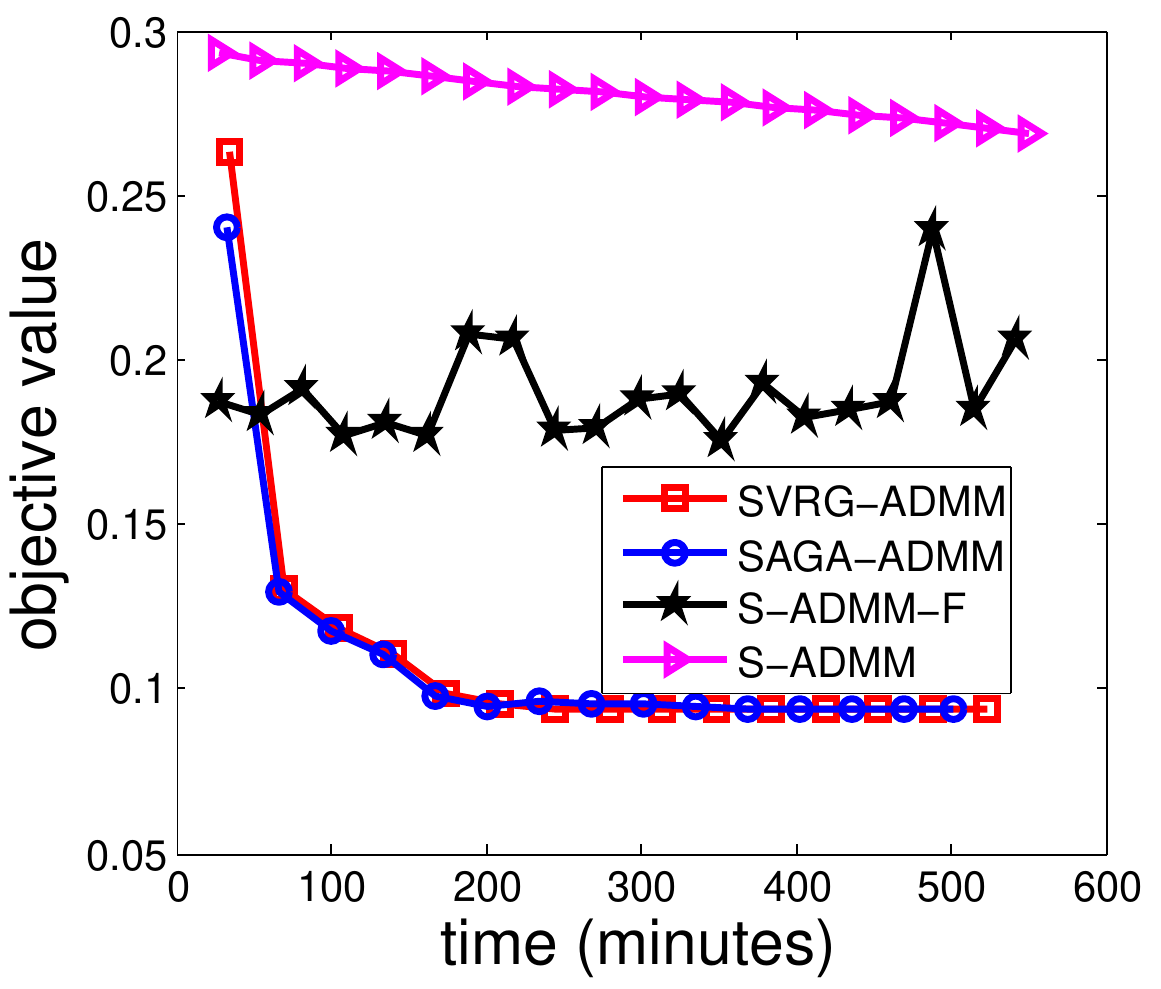}}
\caption{Objective value \emph{versus} time on the \emph{nonconvex} graph-guided model.} \label{1}
\vspace{-1em}
\end{figure*}



\subsection{Experimental Results}

Figures 1-2 show that the test losses of both SVRG-ADMM and SAGA-ADMM faster decrease than that of both S-ADMM and S-ADMM-F,
as the number of effective passes or time increase,
where each effective pass estimates $n$ component gradients.
Figures 3-4 show that the objective values of both SVRG-ADMM and SAGA-ADMM also faster decrease than that of both S-ADMM and S-ADMM-F,
as the number of effective passes or time increase.
In particular, as number of effective passes or time increase, both test loss and objective value of S-ADMM-F are fluctuant.
This implies that the S-ADMM-F may be divergent with a constant $\eta$.
At the same time, though the S-ADMM with time-varying $\eta_t$ is convergence,
but it only slowly converge to the local optimal solution due to
existing of variance of stochastic gradients.
In summary, these experimental results demonstrate the effectiveness of the proposed methods.

\section{Conclusions}

In the paper, we proposed three classes of the nonconvex stochastic ADMM methods with variance reduction,
based on different reduced variance stochastic gradients.
Moreover, we proved that the proposed methods have the iteration complexity bound of $O(1/\epsilon)$
to obtain an $\epsilon$-stationary solution of the nonconvex optimizations.
In particular, we provide a general framework to analyze the convergence
and iteration complexity of stochastic ADMM methods with variance reduction.
Finally, some numerical experiments demonstrate 
the effectiveness of our methods for solving the nonconvex optimizations.
In the future work, we will focus on the study of the convergence 
and iteration complexity of the simple stochastic ADMM method. 
 {\color{blue}{In addition, we will use the momentum acceleration technique
to further accelerate the proposed stochastic ADMM methods.}}


\newpage
\appendix
\section{Proof of Lemma 1}
\label{app:lemma 1}

\begin{proof}
 Since $\hat{\nabla}f(x^{s+1}_t)=\nabla f_{i_t}(x^{s+1}_t)-\nabla f_{i_t}(\tilde{x}^{s})+\nabla f(\tilde{x}^{s})$, we have
 \begin{align}
  \mathbb{E}\|\Delta^{s+1}_t\|^2 & = \mathbb{E}\|\nabla f_{i_t}(x^{s+1}_t)-\nabla f_{i_t}(\tilde{x}^{s})+\nabla f(\tilde{x}^{s})-\nabla f(x_t^{s+1})\|^2 \nonumber \\
  & \mathop{=}^{(i)}  \mathbb{E}\|\nabla f_{i_t}(x^{s+1}_t)-\nabla f_{i_t}(\tilde{x}^{s})\|^2 -
    \|\nabla f(x_t^{s+1})-\nabla f(\tilde{x}^{s})\|^2 \nonumber \\
  & \leq \mathbb{E}\|\nabla f_{i_t}(x^{s+1}_t)-\nabla f_{i_t}(\tilde{x}^{s})\|^2 \nonumber \\
  & = \frac{1}{n}\sum_{i=1}^n \|\nabla f_{i_t}(x^{s+1}_t)-\nabla f_{i_t}(\tilde{x}^{s})\|^2 \nonumber \\
  & \mathop{\leq}^{(ii)} L^2 \|x_t^{s+1}-\tilde{x}^s\|^2. \nonumber
 \end{align}
 where the equality (i) holds by the equality $\mathbb{E}(\xi-\mathbb{E}\xi)^2 = \mathbb{E}\xi^2-(\mathbb{E}\xi)^2$ for random variable $\xi$;
 the inequality (ii) holds by the Assumption 1.
\end{proof}

\section{Proof of the Lemma 2}
\label{app:lemma 2}

\begin{proof}
 Since $\psi_{t}=\frac{1}{n}\sum_{j=1}^n\nabla f_j(z^t_j)$, we have
  \begin{align}
\mathbb{E}\|\Delta_t\|^2 &= \mathbb{E}\|\frac{1}{n} \big(\nabla f_{i_t}(x_t)-\nabla f_{i_t}(z^{t}_{i_t}) \big) -\big(\nabla f(x_t)-\psi_t\big)\|^2 \nonumber \\
  & = \frac{1}{n^2} \mathbb{E}\|\nabla f_{i_t}(x_t)-\nabla f_{i_t}(z^{t}_{i_t})\|^2 + (1-\frac{2}{n})\|\nabla f(x_t)-\psi_t\|^2  \nonumber \\
  & = \frac{1}{n^3} \sum_{i=1}^n \|\nabla f_{i_t}(x_t)-\nabla f_{i_t}(z^{t}_{i})\|^2
  + (1-\frac{2}{n})\frac{1}{n^2}\|\sum_{i=1}^n\big(\nabla f(x_{i_t})-\nabla f_{i_t}(z^{t}_{i})\big)\|^2 \nonumber \\
  & \mathop{\leq}^{(i)} \frac{1}{n^3} \sum_{i=1}^n \|\nabla f_{i_t}(x_t)-\nabla f_{i_t}(z^{t}_{i})\|^2
  + (1-\frac{2}{n})\frac{1}{n}\sum_{i=1}^n\|\big(\nabla f(x_{i_t})-\nabla f_{i_t}(z^{t}_{i})\big)\|^2 \nonumber \\
  & \mathop{\leq}^{(ii)} \big( 1-\frac{1}{n} \big)^2 \frac{L^2}{n} \sum_{i=1}^n \|x_t-z^t_i\|^2. \nonumber
 \end{align}
 where the equality (i) holds by the equality $(\sum_{i=1}^n a_i)^2 \leq n\sum_{i=1}^n a_i^2$;
 the inequality (ii) holds by the Assumption 1.
\end{proof}

\section{Proof of the Lemma 3}
\label{app:lemma 3}

\begin{proof}
 Since $\psi_{t}=\frac{1}{n}\sum_{j=1}^n\nabla f_j(z^t_j)$, we have
  \begin{align}
\mathbb{E}\|\Delta_t\|^2 &= \mathbb{E}\|\nabla f_{i_t}(x_t)-\nabla f_{i_t}(z^{t}_{i_t})+\psi_t-\nabla f(x_t)\|^2 \nonumber \\
  & \mathop{=}^{(i)}  \mathbb{E}\|\nabla f_{i_t}(x_t)-\nabla f_{i_t}(z^{t}_{i_t})\|^2 - \|\nabla f(x_t)-\psi_t\|^2 \nonumber \\
  & \leq \mathbb{E}\|\nabla f_{i_t}(x_t)-\nabla f_{i_t}(z^{t}_{i_t})\|^2 \nonumber \\
  & = \frac{1}{n} \sum_{i=1}^n \|\nabla f_{i_t}(x_t)-\nabla f_{i_t}(z^{t}_{i})\|^2 \nonumber \\
  & \mathop{\leq}^{(ii)} \frac{L^2}{n} \sum_{i=1}^n \|x_t-z^t_i\|^2. \nonumber
 \end{align}
 where the equality (i) holds by the equality $\mathbb{E}(\xi-\mathbb{E}\xi)^2 = \mathbb{E}\xi^2-(\mathbb{E}\xi)^2$ for random variable $\xi$, and
 $\mathbb{E}[\nabla f_{i_t}(z^{t}_{i_t})]=\frac{1}{n}\sum_{j=1}^n\nabla f_j(z^t_j)=\psi_{t}$;
 the inequality (ii) holds by the Assumption 1.
\end{proof}


\section{Proof of Lemma 5}
\label{app:lemma 5}

\begin{proof}
 For notational simplicity, let $x^{s+1}_t = x_t$, $y^{s+1}_t = y_t$, $\lambda^{s+1}_t = \lambda_t$, and $\tilde{x}^s=\tilde{x}$.
 By the optimal condition of step 10 in Algorithm 2, we have
 \begin{align}
  0 & = \hat{\nabla}f(x_t)-A^T\lambda_t+\rho A^T(Ax_{t+1}+By_{t+1}-c)-\eta Q (x_t - x_{t+1}) \nonumber \\
    & = \hat{\nabla}f(x_t)- A^T\lambda_{t+1}-\eta Q (x_t - x_{t+1}), \nonumber
 \end{align}
 where the second equality is due to step 11 in Algorithm 2.
 Thus, we have
 \begin{align}
  A^T\lambda_{t+1} = \hat{\nabla}f(x_t)-\eta Q (x_t - x_{t+1}).
 \end{align}
 By (34), we have
 \begin{align}
 & \|\lambda_{t+1}-\lambda_{t}\|^2 \leq \sigma_{A}^{-1}\|A^T\lambda_{t+1}-A^T\lambda_t\|^2  \nonumber \\
 & \leq \sigma_{A}^{-1}\|\hat{\nabla}f(x_{t})-\hat{\nabla}f(x_{t-1}) - \eta Q (x_t-x_{t+1}) + \eta Q (x_{t-1}-x_t)\|^2 \nonumber \\
 & = \sigma_{A}^{-1}\|\hat{\nabla}f(x_{t})-\nabla f(x_{t}) + \nabla f(x_{t}) - \nabla f(x_{t-1}) + \nabla f(x_{t-1}) - \hat{\nabla}f(x_{t-1})
   -\eta Q (x_t-x_{t+1}) + \eta Q (x_{t-1}-x_t)\|^2 \nonumber \\
 & \mathop{\leq}^{(i)} \frac{5}{\sigma_{A}} \|\hat{\nabla}f(x_{t})-\nabla f(x_{t})\|^2 + \frac{5}{\sigma_{A}} \|\hat{\nabla}f(x_{t-1})-\nabla f(x_{t-1})\|^2
 + \frac{5\eta^2\phi_{\max}^2 }{\sigma_{A}} \|x_t-x_{t+1}\|^2 \nonumber \\
 & \quad + \frac{5(L^2+\eta^2\phi_{\max}^2) }{\sigma_{A}} \|x_{t-1}-x_{t}\|^2,
 \end{align}
 where the inequality (i) holds by the Assumption 1; and $\|Q(x-y)\|^2 \leq \phi^2_{\max}\|x-y\|^2$,
 where $\phi_{\max}$ denotes the largest eigenvalue of positive matrix $Q$.
 Taking expectation conditioned on information $i_t$ to (35), we have
 \begin{align}
\mathbb{E}\|\lambda_{t+1}-\lambda_{t}\|^2 &\leq \frac{5}{\sigma_{A}} \mathbb{E} \|\hat{\nabla}f(x_{t})-\nabla f(x_{t})\|^2
 + \frac{5}{\sigma_{A}} \mathbb{E} \|\hat{\nabla}f(x_{t-1})-\nabla f(x_{t-1})\|^2+ \frac{5\eta^2\phi_{\max}^2}{\sigma_{A}} \|x_t-x_{t+1}\|^2 \nonumber \\
 & \quad + \frac{5(L^2+\eta^2\phi_{\max}^2)}{\sigma_{A}} \|x_{t-1}-x_{t}\|^2 \nonumber \\
 & \mathop{\leq}^{(i)} \frac{5L^2}{\sigma_{A}} \mathbb{E}\|x_{t}-\tilde{x}\|^2 + \frac{5L^2}{\sigma_{A}}\|x_{t-1}-\tilde{x}\|^2
 + \frac{5\eta^2\phi_{\max}^2}{\sigma_{A}} \mathbb{E}\|x_t-x_{t+1}\|^2 + \frac{5(L^2+\eta^2\phi_{\max}^2)}{\sigma_{A}}\|x_{t-1}-x_{t}\|^2, \nonumber
 \end{align}
 where the inequality $(i)$ holds by the Lemma 1.
\end{proof}

\section{Proof of Lemma 6}
\label{app:lemma 6}

\begin{proof}
This proof includes two parts: First, we will prove that the sequence $\{(\Psi^{s}_{t})_{t=1}^m\}_{s=1}^S$
is monotonically decreases over $t\in \{1,2,\cdots,m\}$ in each epoch $s\in \{1,2,\cdots,S\}$;
Second, we will prove that $\Psi^{s}_{m} \geq \Psi^{s+1}_{1}$ for any $s\in \{1,2,\cdots,S\}$.

For notational simplicity, we omit the label $s$  in the first part, i.e.,
let $x^{s+1}_t = x_t$, $y^{s+1}_t = y_t$, $\lambda^{s+1}_t = \lambda_t$, and $\tilde{x}^{s}=\tilde{x}$.
By the step 8 of Algorithm 2, we have
\begin{align}
 \mathcal {L}_\rho (x_t, y_{t+1},\lambda_t) \leq \mathcal {L}_\rho (x_t, y_{t},\lambda_t).
\end{align}
By the optimal condition of step 10 in Algorithm 2, we have
\begin{align}
0 &=(x_t-x_{t+1})^T\big[\hat{\nabla}f(x_t)-A^T\lambda_t+\rho(Ax_{t+1}+By_{t+1}-c) - \eta Q (x_t-x_{t+1})\big] \nonumber \\
 & = (x_t-x_{t+1})^T\big[\hat{\nabla}f(x_t) - \nabla f(x_t) + \nabla f(x_t) -A^T\lambda_t +\rho A^T(Ax_{t+1}+By_{t+1}-c)-\eta Q(x_t-x_{t+1})\big] \nonumber \\
 & \mathop{\leq}^{(i)} f(x_t) - f(x_{t+1}) + (x_t -x_{t+1})^T(\hat{\nabla}f(x_t) - \nabla f(x_t))
  + \frac{L}{2} \|x_{t+1}-x_t\|^2  - \eta\|x_{t+1}-x_t\|^2_Q  \nonumber \\
 & \quad - \lambda_t^T(A x_t-Ax_{t+1}) + \rho (Ax_t -Ax_{t+1})^T(Ax_{t+1}+By_{t+1}-c) \nonumber \\
 & \mathop{=}^{(ii)} f(x_t) - f(x_{t+1}) + (x_t -x_{t+1})^T(\hat{\nabla}f(x_t) - \nabla f(x_t))+\frac{L}{2} \|x_{t+1}-x_t\|^2  - \eta\|x_{t+1}-x_t\|^2_Q \nonumber \\
 & \quad - \lambda_t^T(A x_t+By_{t+1}-c) + \lambda_t^T(Ax_{t+1}+By_{t+1}-c) + \frac{\rho}{2}\|Ax_{t}+By_{t+1}-c\|^2 \nonumber \\
 & \quad - \frac{\rho}{2}\|Ax_{t+1}+By_{t+1}-c\|^2 -\frac{\rho}{2}\|Ax_t -Ax_{t+1}\|^2 \nonumber \\
 & = \mathcal {L}_\rho (x_t, y_{t+1},\lambda_t)- \mathcal {L}_\rho (x_{t+1}, y_{t+1},\lambda_t)
  + (x_t -x_{t+1})^T(\hat{\nabla}f(x_t) - \nabla f(x_t)) \nonumber \\
 & \quad +\frac{L}{2} \|x_{t+1}-x_t\|^2  - \eta\|x_{t+1}-x_t\|^2_Q -\frac{\rho}{2}\|Ax_t -Ax_{t+1}\|^2 \nonumber \\
 & \mathop{\leq}^{(iii)} \mathcal {L}_\rho (x_t, y_{t+1},\lambda_t)- \mathcal {L}_\rho (x_{t+1}, y_{t+1},\lambda_t)
  + (x_t -x_{t+1})^T(\hat{\nabla}f(x_t) - \nabla f(x_t)) \nonumber \\
 & \quad - (\eta\phi_{\min} + \frac{\sigma_A\rho}{2}- \frac{L}{2})\|x_t-x_{t+1}\|^2,
\end{align}
where the inequality (i) holds by the Assumption 1; the equality (ii) holds by using the equality
$(a-b)^T(b-c) = \frac{1}{2}(\|a-c\|^2-\|a-b\|^2-\|b-c\|^2)$ on the term $\rho (Ax_t -Ax_{t+1})^T(Ax_{t+1}+By_{t+1}-c)$;
the inequality (iii) holds by using $-\phi_{\min}\|x_{t+1}-x_t\|^2 \geq -\|x_{t+1}-x_t\|^2_Q$ and $-\sigma_A\|x_{t+1}-x_t\|^2 \geq -\|Ax_t -Ax_{t+1}\|^2$.
Then taking expectation conditioned on information $i_t$ to (37), and using $\mathbb{E}[\hat{\nabla} f(x_{t})]=\nabla f(x_{t})$, we have
\begin{align}
\mathbb{E} [\mathcal {L}_\rho (x_{t+1}, y_{t+1},\lambda_t)] & \leq \mathcal {L}_\rho (x_t, y_{t+1},\lambda_t)
 - (\eta\phi_{\min} + \frac{\sigma_A\rho}{2} - \frac{L}{2})\mathbb{E}\|x_t-x_{t+1}\|^2.
\end{align}

By the step 11 of Algorithm 2, and taking expectation conditioned on information $i_t$, we have
\begin{align}
\mathbb{E} [\mathcal {L}_\rho (x_{t+1}, y_{t+1},\lambda_{t+1})-\mathcal {L}_\rho (x_{t+1}, y_{t+1},\lambda_t)]
&=\frac{1}{\rho}\mathbb{E} \|\lambda_{t+1}-\lambda_t\|^2 \nonumber \\
& \mathop{\leq}^{(i)} \frac{5L^2}{\sigma_{A} \rho} \mathbb{E}\|x_{t}-\tilde{x}\|^2 + \frac{5L^2}{\sigma_{A} \rho}\|x_{t-1}-\tilde{x}\|^2
 + \frac{5\eta^2\phi^2_{\max}}{\sigma_{A} \rho} \mathbb{E}\|x_{t+1}-x_t\|^2 \nonumber \\
& \quad + \frac{5(L^2+\eta^2\phi^2_{\max})}{\sigma_{A} \rho}\|x_{t}-x_{t-1}\|^2,
\end{align}
where the inequality $(i)$ holds by the Lemma 5.

Combining (36), (38) with (39), we have
\begin{align}
\mathbb{E} [\mathcal {L}_\rho (x_{t+1}, y_{t+1},\lambda_{t+1})] \leq & \mathcal {L}_\rho (x_t, y_{t},\lambda_t) + \frac{5L^2}{\sigma_{A}\rho} \mathbb{E}\|x_{t}-\tilde{x}\|^2 + \frac{5L^2}{\sigma_{A}\rho}\|x_{t-1}-\tilde{x}\|^2 + \frac{5(L^2+\eta^2\phi^2_{\max} )}{\sigma_{A}\rho}\|x_{t}-x_{t-1}\|^2 \nonumber \\
 & - (\eta\phi_{\min} + \frac{\sigma_A\rho}{2}- \frac{L}{2} -\frac{5\eta^2 \phi^2_{\max}}{\sigma_{A}\rho})\|x_{t+1}-x_t\|^2.
\end{align}
Next, considering $\mathbb{E}\|x_{t+1}-\tilde{x}\|^2$, we have
\begin{align}
\mathbb{E}\|x_{t+1}-\tilde{x}\|^2 & = \mathbb{E}\|x_{t+1}-x_t+x_t-\tilde{x}\|^2 \nonumber \\
& = \mathbb{E}[\|x_{t+1}-x_t\|^2+ 2(x_{t+1}-x_t)^T(x_t-\tilde{x} )+\|x_t-\tilde{x}\|^2] \nonumber \\
& \leq \mathbb{E}[\|x_{t+1}-x_t\|^2+ 2(\frac{1}{2\beta}\|x_{t+1}-x_t\|^2
+ \frac{\beta }{2}\|x_t-\tilde{x}\|^2 )+\|x_t-\tilde{x}\|^2] \nonumber \\
& = (1+\frac{1}{\beta})\|x_{t+1}-x_t\|^2 + (1+\beta)\|x_t-\tilde{x}\|^2,
\end{align}
where $\beta>0$, and the inequality is due to the Cauchy inequality.
Combining the inequalities (40) and (41), we have
\begin{align}
 &\mathbb{E} [\mathcal {L}_\rho (x_{t+1}, y_{t+1},\lambda_{t+1})+\frac{5(L^2+\eta^2\phi^2_{\max})}{\sigma_{A} \rho}\|x_{t+1}-x_t\|^2
  + h^s_{t+1}(\|x_{t+1}-\tilde{x}\|^2 + \|x_{t}-\tilde{x}\|^2)] \nonumber \\
 & \leq \mathcal {L}_\rho (x_{t}, y_{t},\lambda_{t}) + \frac{5(L^2+\eta^2\phi^2_{\max})}{\sigma_{A} \rho} \|x_t - x_{t-1}\|^2
  + [(2+\beta)h^s_{t+1} + \frac{5L^2}{\sigma_{A} \rho}](\|x_t-\tilde{x}\|^2+\|x_{t-1}-\tilde{x}\|^2)  \nonumber \\
 & \quad -\big[ \eta\phi_{\min} + \frac{\sigma_A\rho}{2}- \frac{L}{2} -\frac{5(2\eta^2\phi^2_{\max}+L^2)}{\sigma_{A}\rho}-(1+\frac{1}{\beta})h_{t+1} \big]\mathbb{E}\|x_{t+1}-x_t\|^2
  - (2+\beta)h^s_{t+1}\|x_{t-1}-\tilde{x}\|^2,
\end{align}
where $h_{t+1}>0$. Then using the definition of the sequence $\{(\Psi_t^s)_{t=1}^m\}_{s=1}^S$, (16) and (17) to (42), we have
\begin{align}
  \Psi^{s+1}_{t+1} \leq \Psi^{s+1}_{t}-\Gamma^{s+1}_t \mathbb{E}\|x^{s+1}_{t+1}-x^{s+1}_t\|^2 - (2+\beta)h^{s+1}_{t+1}\|x^{s+1}_{t-1}-\tilde{x}^{s}\|^2,
\end{align}
for any $s\in \{0,1,\cdots,S-1\}$.
Since $\Gamma^s_t >0, \ \forall t\in \{1,2,\cdots, m\}$, we prove the first part.

Next, we will prove the second part. We begin with considering the upper bound of $\mathbb{E}\|\lambda^{s+1}_0-\lambda^{s+1}_1\|^2$.
Since $\lambda^{s+1}_0 = \lambda^s_m$ and $x^{s+1}_0=x^s_m=\tilde{x}^s$, we have
\begin{align}
 \mathbb{E}\|\lambda^{s+1}_0-\lambda^{s+1}_1\|^2 & = \mathbb{E}\|\lambda^{s}_m-\lambda^{s+1}_1\|^2 \nonumber \\
 & \leq \frac{1}{\sigma_{A}} \mathbb{E}\|A^T\lambda^{s}_m-A^T\lambda^{s+1}_1\|^2 \nonumber \\
 & \mathop{=}^{(i)} \frac{1}{\sigma_{A}} \mathbb{E}\|\hat{\nabla}f(x^s_{m-1})-\hat{\nabla}f(x^{s+1}_{0}) -\eta Q(x^s_{m-1}-x^s_m)
  -\eta Q(x^{s+1}_{0}-x^{s+1}_1)\|^2  \nonumber \\
 & \mathop{=}^{(ii)} \frac{1}{\sigma_{A}} \mathbb{E}\|\hat{\nabla}f(x^s_{m-1})
   - \nabla f(x^s_{m-1}) + \nabla f(x^s_{m-1})-\nabla f(x^{s}_{m})\nonumber \\
 & \quad -\eta Q (x^s_{m-1}-x^s_m) - \eta Q(x^{s+1}_{0}-x^{s+1}_1)\|^2  \nonumber \\
 & \leq \frac{5L^2}{\sigma_{A}} \|x^s_{m-1}-\tilde{x}^{s-1}\|^2 + \frac{5(L^2+\eta^2\phi^2_{\max})}{\sigma_{A}} \|x^s_{m-1}-x^s_{m}\|^2 \nonumber \\
 & \quad + \frac{5\eta^2\phi^2_{\max}}{\sigma_{A}} \|x^{s+1}_{0}-x^{s+1}_1\|^2,
\end{align}
where the equality $(i)$ holds by the equality (34), and the equality $(ii)$ holds by the following result:
\begin{align}
  \hat{\nabla}f(x^{s+1}_{0}) & = \nabla f_{i_t}(x^{s+1}_{0}) - \nabla f_{i_t}(\tilde{x}^{s}) + \nabla f(\tilde{x}^{s}) \nonumber \\
  & = \nabla f_{i_t}(x^{s}_{m}) - \nabla f_{i_t}(x^{s}_m) + \nabla f(x^{s}_m) \nonumber \\
  & = \nabla f(x^{s}_m). \nonumber
\end{align}
By (36), we have
\begin{align}
 \mathcal {L}_\rho (x^{s+1}_0, y^{s+1}_{1},\lambda^{s+1}_0)&\leq \mathcal {L}_\rho (x^{s+1}_0, y^{s+1}_{0},\lambda^{s+1}_0)
  = \mathcal {L}_\rho (x^{s}_m, y^{s}_{m},\lambda^{s}_m).
\end{align}
Similarly, by (38), we have
\begin{align}
\mathbb{E} [\mathcal {L}_\rho (x^{s+1}_{1}, y^{s+1}_{1},\lambda^{s+1}_0)] &\leq \mathcal {L}_\rho (x^{s+1}_0, y^{s+1}_{1},\lambda^{s+1}_0)
 - (\eta\phi_{\min} + \frac{\sigma_A\rho}{2} - \frac{L}{2})\mathbb{E}\| x^{s+1}_{1}-x^{s+1}_0\|^2.
\end{align}
Using (44), we have
\begin{align}
  \mathbb{E} [\mathcal {L}_\rho (x^{s+1}_{1}, y^{s+1}_{1},\lambda^{s+1}_{1})-\mathcal {L}_\rho (x^{s+1}_{1}, y^{s+1}_{1},\lambda^{s+1}_0)]
 & = \frac{1}{\rho} \mathbb{E} \|\lambda^{s+1}_0 - \lambda^{s+1}_1\|^2  \nonumber \\
 & \leq \frac{5L^2}{\sigma_{A}\rho} \|x^s_{m-1}-\tilde{x}^{s-1}\|^2 + \frac{5(L^2+\eta^2\phi^2_{\max})}{\sigma_{A}\rho} \|x^s_{m-1}-x^s_{m}\|^2 \nonumber \\
 & \quad + \frac{5\eta^2\phi^2_{\max}}{\sigma_{A}\rho} \|x^{s+1}_1 - x^{s+1}_{0}\|^2.
\end{align}
Since $\mathcal {L}_\rho (x^{s+1}_0, y^{s+1}_{0},\lambda^{s+1}_0)=\mathcal {L}_\rho (x^{s}_m, y^{s}_{m},\lambda^{s}_m)$,
by combining (45), (46) with (47), we have
\begin{align}
\mathbb{E} [\mathcal {L}_\rho (x^{s+1}_{1}, y^{s+1}_{1},\lambda^{s+1}_{1})] &\leq \mathcal {L}_\rho (x^{s}_m, y^{s}_{m},\lambda^{s}_m)
  + \frac{5L^2}{\sigma_{A}\rho} \|x^s_{m-1}-\tilde{x}^{s-1}\|^2 + \frac{5(L^2+\eta^2\phi^2_{\max})}{\sigma_{A}\rho} \|x^s_{m-1}-x^s_{m}\|^2 \nonumber \\
 & - (\eta\phi_{\min} + \frac{\sigma_A\rho}{2} - \frac{L}{2}-\frac{5\eta^2\phi^2_{\max}}{\sigma_{A}\rho})\mathbb{E}\|x^{s+1}_{1} - x^{s+1}_0\|^2.
\end{align}
Then using $h_1^{s+1}=\frac{10L^2}{\sigma_{A}\rho}$, we have
\begin{align}
 &\mathbb{E} \big[\mathcal {L}_\rho (x^{s+1}_{1}, y^{s+1}_{1},\lambda^{s+1}_{1})
 + h^{s+1}_1[\|x^{s+1}_1-\tilde{x}^s\|^2 + \|x^{s+1}_0-\tilde{x}^s\|^2]
 + \frac{5(L^2+\eta^2\phi^2_{\max})}{\sigma_{A}\rho}\|x^{s+1}_{1}-x^{s+1}_{0}\|^2 \big] \nonumber \\
 & \leq \mathcal {L}_\rho (x^{s}_m, y^{s}_{m},\lambda^{s}_m) + \frac{10L^2}{\sigma_{A}\rho} [\|x^s_{m}-\tilde{x}^{s-1}\|^2 + \|x^s_{m-1}-\tilde{x}^{s-1}\|^2 ]
  + \frac{5(L^2+\eta^2\phi^2_{\max})}{\sigma_{A}\rho} \|x^s_{m-1}-x^s_{m}\|^2 \nonumber \\
 & \quad - (\eta\phi_{\min} + \frac{\sigma_A\rho}{2} - h_1^{s+1}-\frac{L}{2}-\frac{10\eta^2\phi^2_{\max}+5L^2}{\sigma_{A}\rho})\mathbb{E}\|x^{s+1}_0 - x^{s+1}_{1}\|^2
   - \frac{5L^2}{\sigma_{A}\rho} \|x^s_{m-1}-\tilde{x}^{s-1}\|^2  \nonumber \\
 & \quad -\frac{10L^2}{\sigma_{A}\rho} \|x^s_{m}-\tilde{x}^{s-1}\|^2.
\end{align}
Finally, using the definition of the sequence $\{(\Psi_t^s)_{t=1}^m\}_{s=1}^S$, (16) and (17) to (49), we have
\begin{align}
 \Psi^{s+1}_{1} \leq \Psi^{s}_{m} - \Gamma^{s}_m\mathbb{E}\|x^{s+1}_0 - x^{s+1}_{1}\|^2 -\frac{5L^2}{\sigma_{A}\rho} [\|x^s_{m-1}-\tilde{x}^{s-1}\|^2.
\end{align}
Since $\Gamma^{s}_m>0,\ \forall s \geq 1$, we can obtain the above result of the second part.
Thus, we prove the above conclusion.
\end{proof}

\section{Proof of Lemma 7}
\label{app:lemma 7}

\begin{proof}
 By definition of the sequence $\{(\Psi_t^s)_{t=1}^m\}_{s=1}^S$, we have
 \begin{align}
   \Psi^{s}_{t}& \geq \mathbb{E}[\mathcal {L}_\rho (x^{s}_{t}, y^{s}_{t},\lambda^{s}_{t})] \nonumber \\
  & = f(x^{s}_{t})+g(y^{s}_{t}) - (\lambda^{s}_{t})^T(Ax^{s}_{t}+By^{s}_{t}-c)
  + \frac{\rho}{2}\|Ax^{s}_{t}+By^{s}_{t}-c\|^2 \nonumber \\
  & \mathop{=}^{(i)} f(x^{s}_{t})+g(y^{s}_{t}) - \frac{1}{\rho}(\lambda^{s}_{t})^T(\lambda^{s}_{t-1} -\lambda^{s}_{t})
   + \frac{1}{2\rho}\|\lambda^{s}_{t-1} - \lambda^{s}_{t}\|^2 \nonumber \\
  & = f(x^{s}_{t})+g(y^{s}_{t}) - \frac{1}{2\rho}\|\lambda^{s}_{t-1}\|^2 + \frac{1}{2\rho}\|\lambda^{s}_{t}\|^2
  + \frac{1}{\rho}\|\lambda_t-\lambda_{t-1}\|^2\nonumber \\
  & \mathop{\geq}^{(ii)} f^* + g^* - \frac{1}{2\rho}\|\lambda^{s}_{t-1}\|^2 + \frac{1}{2\rho}\|\lambda^{s}_{t}\|^2,
 \end{align}
 where the equality $(i)$ holds by the step 11 of Algorithm 2,
 and the inequality $(ii)$ holds by Assumption 2.

 Summing the inequality (51) over $t=1,2,\cdots,m$ and $s=1,2,\cdots,S$, we have
 \begin{align}
 & \frac{1}{T}\sum_{s=1}^{S}\sum_{t=1}^{m} \Psi^{s}_{t} \geq f^* + g^* - \frac{1}{2\rho}\|\lambda^1_{0}\|^2. \nonumber
 \end{align}
 Therefore, we can obtain the above result.
\end{proof}

\section{Proof of Theorem 8}
\label{app:theorem 8}

\begin{proof}
Using the above inequalities (43) and (50), we have
\begin{align}
 \Psi_{t+1}^{s}  \leq \Psi_{t}^{s}-\Gamma_t^s\mathbb{E}\|x_{t+1}^{s}-x_t^{s}\|^2
 - (2+\beta)h^{s}_{t+1}\|x^{s}_{t-1}-\tilde{x}^{s-1}\|^2,
\end{align}
and
\begin{align}
\Psi^{s+1}_{1}  \leq \Psi^{s}_{m} - \Gamma^{s}_m\mathbb{E}\|x^{s+1}_0 - x^{s+1}_{1}\|^2
 -\frac{5L^2}{\sigma_{A}\rho} \|x^s_{m}-\tilde{x}^{s-1}\|^2.
\end{align}
for any $s\in \{1,2,\cdots, S\}$ and $t\in\{1,2,\cdots,m\}$.
Summing (52) and (53) over $t=1,2,\cdots,m$ and $s=1,2,\cdots,S$, we have
\begin{align}
 \Psi_{m}^{S}- \Psi_{1}^{1}
 \leq -\gamma\sum_{s=1}^{S}\sum_{t=1}^{m}\mathbb{E}\|x_{t}^{s}-x_{t-1}^{s}\|^2
  - \omega \sum_{s=1}^{S}\sum_{t=1}^{m} \|x_{t-1}^{s}-\tilde{x}^{s-1}\|^2
\end{align}
where $\gamma = \min_{s,t} \Gamma^s_t$,
and $\omega = \min_{s,t}\{(2+\beta)h^{s}_{t+1},\frac{5L^2}{\sigma_{A} \rho} \} = \frac{5L^2}{\sigma_{A} \rho}$.
From Lemma 7, there exists a low bound $\Psi^*$ of the sequence $\{\Psi^s_t\}$, i.e., $ \Psi^{s}_{t} \geq \Psi^*$.
Using (54) and the definition of $\theta^s_t$, we have
\begin{align}
 \theta^{\hat{s}}_{\hat{t}} = \min_{s,t} \theta^{s}_t \leq \frac{1}{\tau T} (\Psi^{1}_{1} - \Psi^*),
\end{align}
where $\tau = \min(\gamma,\omega)$, and $T=mS$, so $\theta^{\hat{s}}_{\hat{t}} = O(\frac{1}{T})$.

Next, we give the upper bounds to the terms in (11-13) by using $\theta^s_t$.
By (34), we have
\begin{align}
 & \mathbb{E}\|A^T\lambda^{s}_{t}-\nabla f(x^{s}_{t})\|^2  \nonumber \\
 & = \mathbb{E}\|\hat{\nabla}f(x^{s}_{t-1})-\nabla f(x^{s}_{t})-\eta Q(x^{s}_{t-1}-x^{s}_{t})\|^2 \nonumber \\
 & = \mathbb{E}\|\hat{\nabla}f(x^{s}_{t-1})-\nabla f(x^{s}_{t-1}) +\nabla f(x^{s}_{t-1})- \nabla f(x^{s}_{t})
 -\eta Q(x^{s}_{t-1}-x^{s}_{t})\|^2 \nonumber \\
 & \leq 3L^2\|x^{s}_{t-1}-\tilde{x}^{s-1}\|^2  + 3(L^2+\eta^2\phi^2_{\max})\|x^{s}_{t-1}-x^{s}_t\|^2 \nonumber \\
 & \leq 3(L^2+\eta^2\phi^2_{\max})\theta^{s}_t.
\end{align}
By using the step 11 of Algorithm 2 and the Lemma 5, we have
\begin{align}
 \mathbb{E}\|Ax^{s}_{t+1}+By^{s}_{t+1}-c\|^2 & = \frac{1}{\rho^2}\|\lambda^{s}_{t+1}-\lambda^{s}_{t}\|^2  \nonumber \\
 & \leq \frac{5L^2}{\sigma_{A}\rho^2} \mathbb{E} \|x_{t}^{s}-\tilde{x}^{s-1}\|^2
 + \frac{5L^2}{\sigma_{A}\rho^2}\|x_{t-1}^{s}-\tilde{x}^{s-1}\|^2  \nonumber \\
 & + \frac{5\eta^2\phi^2_{\max}}{\sigma_{A}\rho^2} \mathbb{E}\|x_{t+1}^{s}-x_t^{s}\|^2
 + \frac{5(L^2+\eta^2\phi^2_{\max})}{\sigma_{A}\rho^2}\|x_{t}^{s}-x_{t-1}^{s}\|^2 \nonumber \\
 & \leq \frac{5(L^2+\eta^2\phi^2_{\max})}{\sigma_{A}\rho^2} \theta^s_{t}.
\end{align}
By the step 8 of Algorithm 2,
there exists a sub-gradient $\mu \in \partial g(y_{t}^{s})$ such that
\begin{align}
 \mathbb{E}\big[\mbox{dist}( B^T\lambda_{t}^{s}, \partial g(y_{t}^{s}))\big]^2  &\leq \|\mu-B^T\lambda_{t}^{s}\|^2 \nonumber \\
 & = \|B^T\lambda^s_{t-1}-\rho B^T(Ax^s_{t-1}+By_{t}^{s}-c)-B^T\lambda^{s}_{t}\|^2 \nonumber \\
 & = \|\rho B^TA(x^{s}_{t}-x^{s}_{t-1})\|^2 \nonumber \\
 & \leq \rho^2\|B\|_2^2\|A\|_2^2\|x^{s}_{t}-x^{s}_{t-1}\|^2 \nonumber \\
 & \leq \rho^2\|B\|_2^2\|A\|_2^2 \theta^s_{t}.
\end{align}
Finally, using the Definition 4 and (20), we can conclude that the SVRG-ADMM converges an $\epsilon$-stationary point of the problem (1).
\end{proof}


\section{Proof of the Lemma 9}
\label{app:lemma 9}

\begin{proof}
 By the optimal condition of step 7 in Algorithm 3, we have
 \begin{align}
  0 & = \hat{\nabla}f(x_t)-A^T\lambda_t+\rho A^T(Ax_{t+1}+By_{t+1}-c) -\eta Q(x_t - x_{t+1}) \nonumber \\
  & = \hat{\nabla}f(x_t)- A^T\lambda_{t+1}- \eta Q (x_t - x_{t+1}), \nonumber
 \end{align}
 where the second equality is due to step 8 in Algorithm 3.
 Thus, we have
 \begin{align}
 A^T\lambda_{t+1} = \hat{\nabla}f(x_t)- \eta Q (x_t - x_{t+1}).
 \end{align}
 By (59), we have
 \begin{align}
 \|\lambda_{t+1}-\lambda_{t}\|^2 & \leq \sigma_{A}^{-1}\|A^T\lambda_{t+1}-A^T\lambda_t\|^2  \nonumber \\
 & \leq \sigma_{A}^{-1}\|\hat{\nabla}f(x_t)-\hat{\nabla}f(x_{t-1})- \eta Q (x_t-x_{t+1}) + \eta Q (x_{t-1}-x_t)\|^2 \nonumber \\
 & = \sigma_{A}^{-1}\|\hat{\nabla}f(x_{t})-\nabla f(x_{t}) + \nabla f(x_{t}) - \nabla f(x_{t-1}) + \nabla f(x_{t-1})
    - \hat{\nabla}f(x_{t-1}) \nonumber \\
 & \quad -\eta Q (x_t-x_{t+1}) + \eta Q (x_{t-1}-x_t)\|^2 \nonumber \\
 & \mathop{\leq}^{i} \frac{5}{\sigma_{A}} \|\hat{\nabla}f(x_{t})-\nabla f(x_{t})\|^2
 + \frac{5}{\sigma_{A}} \|\hat{\nabla}f(x_{t-1})-\nabla f(x_{t-1})\|^2
 + \frac{5\eta^2\phi^2_{\max}}{\sigma_{A}} \|x_t-x_{t+1}\|^2 \nonumber \\
 & \quad + \frac{5(\eta^2\phi^2_{\max}+L^2)}{\sigma_{A}} \|x_{t-1}-x_{t}\|^2,
 \end{align}
 where the inequality (i) holds by the Assumption 1.

 Taking expectation conditioned on information $i_t$ to (60), we have
 \begin{align}
\mathbb{E}\|\lambda_{t+1}-\lambda_{t}\|^2 &\leq \frac{5}{\sigma_{A}} \mathbb{E} \|\hat{\nabla}f(x_{t})-\nabla f(x_{t})\|^2
  + \frac{5}{\sigma_{A}} \mathbb{E} \|\hat{\nabla}f(x_{t-1})-\nabla f(x_{t-1})\|^2 \nonumber \\
 & \quad + \frac{5\eta^2\phi^2_{\max}}{\sigma_{A}} \|x_t-x_{t+1}\|^2 + \frac{5(\eta^2\phi^2_{\max}+L^2)}{\sigma_{A}} \|x_{t-1}-x_{t}\|^2 \nonumber \\
 & \mathop{\leq}^{(i)} \big(1-\frac{1}{n}\big)^2 \frac{5L^2}{\sigma_{A} n} \sum_{i=1}^n\mathbb{E}\|x_{t}-z^{t}_i\|^2
 + \big(1-\frac{1}{n}\big)^2\frac{5L^2}{\sigma_{A} n}\sum_{i=1}^n \mathbb{E} \|x_{t-1}-z^{t-1}_i\|^2 \nonumber \\
 & \quad  + \frac{5\eta^2\phi^2_{\max}}{\sigma_{A}}\|x_{t+1}-x_t\|^2+ \frac{5(\eta^2\phi^2_{\max} + L^2)}{\sigma_{A}}\|x_{t}-x_{t-1}\|^2, \nonumber
 \end{align}
 where the inequality $(i)$ holds by the Lemma 2.
\end{proof}

\section{Proof of the Lemma 10}
\label{app:lemma 10}

\begin{proof}
By the step 5 of Algorithm 3, we have
\begin{align}
 \mathcal {L}_\rho (x_t, y_{t+1},\lambda_t) \leq \mathcal {L}_\rho (x_t, y_{t},\lambda_t).
\end{align}
Considering the optimal condition of step 7 in Algorithm 3, we have
\begin{align}
0 & =(x_t-x_{t+1})^T\big[\hat{\nabla}f(x_t)+\rho A^T(Ax_{t+1}+By_{t+1}-c) -A^T\lambda_t-\eta Q(x_t-x_{t+1})\big] \nonumber \\
 & = (x_t-x_{t+1})^T\big[\hat{\nabla}f(x_t) - \nabla f(x_t) + \nabla f(x_t) -A^T\lambda_t - \eta Q(x_t-x_{t+1})
 +\rho A^T(Ax_{t+1}+By_{t+1}-c)\big] \nonumber \\
 & \mathop{\leq}^{(i)} f(x_t) - f(x_{t+1}) + (x_t -x_{t+1})^T(\hat{\nabla}f(x_t) - \nabla f(x_t)) + \frac{L}{2}\|x_{t+1}-x_t\|^2 -\eta\|x_{t+1}-x_t\|^2_Q  \nonumber \\
 &\quad - \lambda_t^T(Ax_{t+1}-A x_t) + \rho (Ax_t -Ax_{t+1})^T(Ax_{t+1}+By_{t+1}-c) \nonumber \\
 & \mathop{=}^{(ii)} f(x_t) - f(x_{t+1}) + (x_t -x_{t+1})^T(\hat{\nabla}f(x_t) - \nabla f(x_t)) + \frac{L}{2}\|x_{t+1}-x_t\|^2 -\eta\|x_{t+1}-x_t\|^2_Q \nonumber \\
 & \quad - \lambda_t^T(A x_t+By_{t+1}-c) + \lambda_t^T(Ax_{t+1}+By_{t+1}-c)
  + \frac{\rho}{2}\|Ax_{t}+By_{t+1}-c\|^2 \nonumber \\
 & \quad - \frac{\rho}{2}\|Ax_{t+1}+By_{t+1}-c\|^2 -\frac{\rho}{2}\|Ax_t-Ax_{t+1}\|^2 \nonumber \\
 & = \mathcal {L}_\rho (x_t, y_{t+1},\lambda_t)- \mathcal {L}_\rho (x_{t+1}, y_{t+1},\lambda_t)
 + (x_t -x_{t+1})^T(\hat{\nabla}f(x_t) - \nabla f(x_t)) \nonumber \\
 & \quad + \frac{L}{2}\|x_{t+1}-x_t\|^2 -\eta\|x_{t+1}-x_t\|^2_Q -\frac{\rho}{2}\|Ax_t-Ax_{t+1}\|^2 \nonumber \\
 & \mathop{\leq}^{(iii)} \mathcal {L}_\rho (x_t, y_{t+1},\lambda_t)- \mathcal {L}_\rho (x_{t+1}, y_{t+1},\lambda_t)
 + (x_t -x_{t+1})^T(\hat{\nabla}f(x_t) - \nabla f(x_t)) \nonumber \\
 & \quad - (\eta\phi_{\min} + \frac{\sigma_A\rho}{2}-\frac{L}{2} )\|x_t-x_{t+1}\|^2 \nonumber \\
 & \mathop{\leq}^{(vi)} \mathcal {L}_\rho (x_t, y_{t+1},\lambda_t)- \mathcal {L}_\rho (x_{t+1}, y_{t+1},\lambda_t)
 + \frac{1}{2\vartheta}\|x_t -x_{t+1}\|^2 + \frac{\vartheta}{2}\|\hat{\nabla}f(x_t) - \nabla f(x_t)\|^2\nonumber \\
 & \quad - (\eta\phi_{\min} + \frac{\sigma_A\rho}{2}-\frac{L}{2} )\|x_t-x_{t+1}\|^2
\end{align}
where the inequality (i) holds by the Assumption 1; the equality (ii) holds by using the equality
$(a-b)^T(b-c) = \frac{1}{2}(\|a-c\|^2-\|a-b\|^2-\|b-c\|^2)$ on the term $\rho (Ax_t -Ax_{t+1})^T(Ax_{t+1}+By_{t+1}-c)$;
the inequality (iii) holds by using $-\phi_{\min}\|x_{t+1}-x_t\|^2 \geq -\|x_{t+1}-x_t\|^2_Q$ and $-\sigma_A\|x_{t+1}-x_t\|^2 \geq -\|Ax_t -Ax_{t+1}\|^2$;
the inequality (vi) holds by the Cauchy inequality.

Then taking expectation conditioned on information $i_t$ to (62), we have
\begin{align}
\mathbb{E} [\mathcal {L}_\rho (x_{t+1}, y_{t+1},\lambda_t)] & \leq \mathcal {L}_\rho (x_t, y_{t+1},\lambda_t) - (\eta\phi_{\min} + \frac{\sigma_A\rho}{2} -\frac{L}{2}-\frac{1}{2\vartheta})\|x_t-x_{t+1}\|^2 \nonumber \\
& \quad + \mathbb{E} [\|\hat{\nabla}f(x_t) - \nabla f(x_t)\|^2] \nonumber \\
& \mathop{\leq}^{(i)} \mathcal {L}_\rho (x_t, y_{t+1},\lambda_t)
 - (\eta\phi_{\min} + \frac{\sigma_A\rho}{2} -\frac{L}{2}-\frac{1}{2\vartheta})\|x_t-x_{t+1}\|^2\nonumber \\
&\quad + \big( 1-\frac{1}{n} \big)^2 \frac{\vartheta L^2}{2 n} \sum_{i=1}^n \|x_t-z^t_i\|^2,
\end{align}
where the inequality (i) holds by Lemma 2.
By the step 8 of Algorithm 3, and taking expectation conditioned on information $i_t$, we have
\begin{align}
\mathbb{E}[\mathcal {L}_\rho (x_{t+1}, y_{t+1},\lambda_{t+1})- & \mathcal {L}_\rho (x_{t+1}, y_{t+1},\lambda_t)]
 =  \frac{1}{\rho}\mathbb{E}[\|\lambda_t-\lambda_{t+1}\|^2]\nonumber \\
& \mathop{\leq}^{(i)} \big( 1-\frac{1}{n} \big)^2\frac{5L^2}{\sigma_{A} \rho n} \sum_{i=1}^n \mathbb{E}\|x_{t}-z^{t}_i\|^2
 + \big( 1-\frac{1}{n} \big)^2\frac{5L^2}{\sigma_{A}\rho n}\sum_{i=1}^n \mathbb{E}\|x_{t-1}-z^{t-1}_i\|^2\nonumber \\
& \quad + \frac{5(\eta^2\phi^2_{\max})}{\sigma_{A} \rho}\|x_{t+1}-x_t\|^2 + \frac{5(\eta^2\phi^2_{\max}+L^2)}{\sigma_{A} \rho}\|x_{t}-x_{t-1}\|^2,
\end{align}
where the inequality $(i)$ holds by the Lemma 9.

Combining (61), (63) with (64), we have
\begin{align}
 \mathbb{E} [\mathcal {L}_\rho (x_{t+1}, y_{t+1},\lambda_{t+1})] \leq  & \mathcal {L}_\rho (x_t, y_{t},\lambda_t)
 + \big( 1-\frac{1}{n} \big)^2 \big(\frac{5L^2}{\sigma_{A}\rho n} + \frac{\vartheta L^2}{2 n}\big)\sum_{i=1}^n \mathbb{E}\|x_{t}-z^{t}_i\|^2 \nonumber \\
 & + \big( 1-\frac{1}{n} \big)^2 \frac{5L^2}{\sigma_{A}\rho n}\sum_{i=1}^n \mathbb{E}\|x_{t-1}-z^{t-1}_i\|^2 + \frac{5(\eta^2\phi^2_{\max}+L^2)}{\sigma_{A}\rho}\|x_{t}-x_{t-1}\|^2 \nonumber \\
 & - \big(\eta\phi_{\min} + \frac{\sigma_A\rho}{2} - \frac{L}{2} -\frac{1}{2\vartheta} - \frac{5\eta^2\phi^2_{\max}}{\sigma_{A}\rho} \big)\|x_{t+1}-x_t\|^2.
\end{align}
Next, considering $\frac{1}{n}\sum_{i=1}^n \mathbb{E}\|x_{t+1}-z^{t+1}_i\|^2$, we have
\begin{align}
\frac{1}{n}\sum_{i=1}^n \mathbb{E}\|x_{t+1}-z^{t+1}_i\|^2 =\frac{1}{n}\sum_{i=1}^n \big[ \frac{1}{n}\mathbb{E}\|x_{t+1}-x_{t}\|^2 + \frac{n-1}{n}\mathbb{E}\|x_{t+1} - z^{t}_i\|^2 \big].
\end{align}
The term $\mathbb{E}\|x_{t+1} - z^{t}_i\|^2$ in (66) can be bounded as follows:
\begin{align}
\mathbb{E}\|x_{t+1} - z^{t}_i\|^2 & = \mathbb{E}\|x_{t+1}-x_{t}+x_{t}-z^{t}_i\|^2 \nonumber \\
& = \mathbb{E}[\|x_{t+1}-x_{t}\|^2+ 2(x_{t+1}-x_{t})^T(x_{t}-z^{t}_i ) +\|x_{t-1}-z^{t}_i\|^2] \nonumber \\
& \leq \mathbb{E}[\|x_{t+1}-x_{t}\|^2 + 2(\frac{1}{2\beta}\mathbb{E}\|x_{t+1}-x_{t}\|^2 + \frac{\beta }{2}\|x_{t}-z^{t}_i\|^2 )
 +\|x_{t}-z^{t}_i\|^2] \nonumber \\
& = (1+\frac{1}{\beta})\mathbb{E}\|x_{t+1}-x_{t}\|^2 + (1+\beta)\|x_{t}-z^{t}_i\|^2, \nonumber
\end{align}
where $\beta>0$, and the inequality is due to the Cauchy inequality.
Then, we have
\begin{align}
 &\frac{1}{n}\sum_{i=1}^n \mathbb{E}\|x_{t+1}-z^{t+1}_i\|^2 \leq (1+\frac{(n-1)}{n\beta})\mathbb{E}\|x_{t+1}-x_{t}\|^2  +(1+\beta)\frac{n-1}{n^2}\sum_{i=1}^n \|x_{t}-z^{t}_i\|
\end{align}
Combining (65) with (67), we have
\begin{align}
 &\mathbb{E} \big[\mathcal {L}_\rho (x_{t+1}, y_{t+1},\lambda_{t+1}) + \frac{5(L^2 + \eta^2\phi^2_{\max})}{\sigma_{A} \rho}\|x_{t+1}-x_t\|^2
  + \big( 1-\frac{1}{n} \big)^2 \frac{\alpha_{t+1}}{n}\sum_{i=1}^n (\|x_{t+1}-z^{t+1}_i\|^2 + \|x_{t}-z^{t}_i\|^2) \big ] \nonumber \\
 & \leq \mathcal {L}_\rho (x_{t}, y_{t},\lambda_{t}) + \frac{5(L^2+\eta^2\phi^2_{\max})}{\sigma_{A} \rho} \|x_t - x_{t-1}\|^2 \nonumber\\
 & \quad + \big( 1-\frac{1}{n} \big)^2 \big[(2+\beta-\frac{1+\beta}{n})\alpha_{t+1} + \frac{5L^2}{\sigma_{A} \rho}
  + \frac{\vartheta L^2}{2}\big]\frac{1}{n}\sum_{i=1}^n \big( \|x_{t}-z^{t}_i\|^2
  + \|x_{t-1}-z^{t-1}_i\|^2 \big) \nonumber\\
 & \quad -\big[\eta\phi_{\min} + \frac{\sigma_A\rho}{2}- \frac{L}{2} - \frac{1}{2\vartheta}-\frac{5(2\eta^2\phi^2_{\max} + L^2)}{\sigma_{A} \rho}
   -\big( 1-\frac{1}{n} \big)^2(1+\frac{1}{\beta}-\frac{1}{n\beta})\alpha_{t+1}\big]\|x_{t+1}-x_t\|^2 \nonumber\\
 & \quad -\big( 1-\frac{1}{n} \big)^2\big[(2+\beta-\frac{1+\beta}{n})\alpha_{t+1} +  \frac{\vartheta L^2}{2} \big ]\frac{1}{n}\sum_{i=1}^n \|x_{t-1}-z^{t-1}_i\|^2.
\end{align}

Finally, using the definition of the sequence $\{\Phi_t\}_{t=1}^T$, (22) and (23), we have
\begin{align}
\Phi_{t+1} & \leq  \Phi_{t}-\Gamma_{t}\|x_{t+1}-x_t\|^2
  -\big( 1-\frac{1}{n} \big)^2\big[(2+\beta-\frac{1+\beta}{n})\alpha_{t+1} +  \frac{\vartheta L^2}{2} \big ]\frac{1}{n}\|x_{t-1}-z^{t-1}_i\|^2.
\end{align}
Since $\Gamma_t >0$ for any $t\in\{1,2,\cdots,T\}$, we can obtain the above result.
\end{proof}

\section{Proof of the Lemma 11}
\label{app:lemma 11}

\begin{proof}
 By the definition of the sequence $\{\Phi_t\}_{t=1}^T$, we have
 \begin{align}
  \Phi(x_{t},y_{t},\lambda_{t},z^{t})
  &\geq \mathbb{E}[\mathcal {L}_\rho (x_{t}, y_{t},\lambda_{t})] \nonumber \\
  & = f(x_{t})+g(y_{t}) - \lambda_{t}^T(Ax_{t}+By_{t}-c) + \frac{\rho}{2}\|Ax_{t}+By_{t}-c\|^2 \nonumber \\
  & \mathop{=}^{(i)} f(x_{t})+g(y_{t}) - \frac{1}{\rho}\lambda_{t}^T(\lambda_{t-1} -\lambda_{t})
  + \frac{1}{2\rho}\|\lambda_{t-1}-\lambda_{t}\|^2  \nonumber \\
  & = f(x_{t})+g(y_{t}) - \frac{1}{2\rho}\|\lambda_{t-1}\|^2 + \frac{1}{2\rho}\|\lambda_{t}\|^2
  + \frac{1}{\rho}\|\lambda_t-\lambda_{t-1}\|^2\nonumber \\
  & \mathop{\geq}^{(ii)} f^* + g^* - \frac{1}{2\rho}\|\lambda_{t-1}\|^2 + \frac{1}{2\rho}\|\lambda_{t}\|^2,
 \end{align}
 where the equality $(i)$ holds by the step 8 in Algorithm 3, and the inequality $(ii)$ holds by Assumption 2.

 Summing the inequality (70) over $t=1,1,\cdots,T$, we have
 \begin{align}
 & \frac{1}{T}\sum_{t=1}^{T} \Phi_{t} \geq f^* + g^* - \frac{1}{2\rho}\|\lambda_{0}\|^2. \nonumber
 \end{align}
 Therefore, we can obtain the above result.
\end{proof}

\section{Proof of the Theorem 12}
\label{app:theorem 12}

\begin{proof}
Using (69), we have
\begin{align}
 \Phi_{t+1}  \leq \Phi_{t}-\Gamma_{t+1}\|x_{t+1}-x_t\|^2 - \big( 1-\frac{1}{n} \big)^2\big[(2+\beta-\frac{1+\beta}{n})\alpha_{t+1}
  +  \frac{\vartheta L^2}{2} \big ]\frac{1}{n}\|x_{t-1}-z^{t-1}_i\|^2,
\end{align}
for $t\in\{1,2,\cdots,T\}$.
Summing the inequality (71) over $t=1,2,\cdots,T$, we have
\begin{align}
 \Phi_{T} &\leq \Phi_{1} -\gamma\sum_{t=1}^{T} \mathbb{E}\|x_{t+1}-x_t\|^2
   - \omega \sum_{t=1}^T\frac{1}{n}\sum_{i=1}^n\|x_{t-1}-z^{t-1}_i\|^2.
\end{align}
where $\gamma = \min_t \Gamma_t $ and $\omega=\min_t \big( 1-\frac{1}{n} \big)^2\big[(2+\beta-\frac{1+\beta}{n})\alpha_{t+1} + \frac{\vartheta L^2}{2}\big]$.
From Lemma 11, there exists a low bound $\Phi^*$ of the sequence $\{\Phi_t\}$, such that $ \Phi_{t}\geq \Phi^*$ for $\forall t \geq 1$.
Then by the definition of $\theta_{t}$ and (72), we have
\begin{align}
\theta_{\hat{t}}=\min_{1 \leq t \leq T} \theta_{t} \leq \frac{1}{\tau T} (\Phi_{1}- \Phi^*),
\end{align}
where $\tau=\min(\gamma,\omega)$ , so $\theta_{\hat{t}} = O(\frac{1}{T})$.

Next, we give upper bounds to the terms in (11-13) by using $\theta_{t}$.
By (59), we have
\begin{align}
 & \mathbb{E}\|A^T\lambda_{t+1}-\nabla f(x_{t+1})\|^2  \nonumber \\
 & = \mathbb{E}\|\hat{\nabla}f(x_{t}) - \nabla f(x_{t+1}) - \eta Q (x_t-x_{t+1})\|^2 \nonumber \\
 & = \mathbb{E}\|\hat{\nabla}f(x_{t})-\nabla f(x_{t}) +\nabla f(x_{t})- \nabla f(x_{t+1})
  - \eta Q (x_t-x_{t+1})\|^2  \nonumber \\
 & \leq \big(1-\frac{1}{n}\big)^2 \frac{3L^2}{n}\sum_{i=1}^n\|x_{t}-z^t_i\|^2  + 3(L^2+\eta^2\phi^2_{\max})\|x_t-x_{t+1}\|^2 \nonumber \\
 & \leq 3(L^2+\eta^2\phi_{\max}^2)\theta_{t}.
\end{align}
By the step 8 of Algorithm 3 and the Lemma 8, we have
\begin{align}
\mathbb{E}\|Ax_{t+1}+By_{t+1}-c\|^2 &= \frac{1}{\rho^2}\|\lambda_{t+1}-\lambda_t\|^2  \nonumber \\
 & \leq \big(1-\frac{1}{n}\big)^2\frac{5L^2}{\sigma_{A}\rho^2 n} \sum_{i=1}^n \mathbb{E} \|x_{t+1}-z^{t+1}_i\|^2
 + \big(1-\frac{1}{n}\big)^2 \frac{5L^2}{\sigma_{A}\rho^2 n} \sum_{i=1}^n \|x_{t}-z^t_i\|^2  \nonumber \\
 & \quad + \frac{5L^2}{\sigma_{A}\rho^2} \mathbb{E} \|x_{t+1}-x_t\|^2
  + \frac{5(L^2+\eta^2\phi^2_{\max})}{\sigma_{A}\rho^2}\|x_{t}-x_{t-1}\|^2  \nonumber \\
 & \leq \frac{5(L^2+\eta^2\phi_{\max}^2)}{\sigma_{A}\rho^2} \theta_{t}.
\end{align}
By the step 5 of Algorithm 3,
there exists a subgradient $\mu \in \partial g(y_{t+1})$ such that
\begin{align}
 \mathbb{E}\big[\mbox{dist} (B^T\lambda_{t+1}, \partial g(y_{t+1}))\big]^2  & \leq \|\mu-B^T\lambda_{t+1}\|^2 \nonumber \\
 & = \|B^T\lambda_t-\rho B^T(Ax_t+By_{t+1}-c)-B^T\lambda_{t+1}\|^2 \nonumber \\
 & = \|\rho B^TA(x_{t+1}-x_{t})\|^2 \nonumber \\
 & \leq \rho^2\|B\|_2^2\|A\|_2^2\|x_{t+1}-x_t\|^2 \nonumber \\
 & \leq \rho^2\|B\|_2^2\|A\|_2^2 \theta_{t}.
\end{align}
Thus, by (26) and the Definition 4, we conclude that the SAG-ADMM can converge an $\epsilon$-stationary point of
the problem (1).
\end{proof}


\section{Proof of the Lemma 13}
\label{app:lemma 13}

\begin{proof}
 By the optimal condition of step 7 in Algorithm 4, we have
 \begin{align}
  0 & = \hat{\nabla}f(x_t)-A^T\lambda_t+\rho A^T(Ax_{t+1}+By_{t+1}-c) -\eta Q(x_t - x_{t+1}) \nonumber \\
  & = \hat{\nabla}f(x_t)- A^T\lambda_{t+1}- \eta Q (x_t - x_{t+1}), \nonumber
 \end{align}
 where the second equality is due to step 8 in Algorithm 4.
 Thus, we have
 \begin{align}
 A^T\lambda_{t+1} = \hat{\nabla}f(x_t)- \eta Q (x_t - x_{t+1}).
 \end{align}
 By (77), we have
 \begin{align}
 \|\lambda_{t+1}-\lambda_{t}\|^2 & \leq \sigma_{A}^{-1}\|A^T\lambda_{t+1}-A^T\lambda_t\|^2  \nonumber \\
 & \leq \sigma_{A}^{-1}\|\hat{\nabla}f(x_t)-\hat{\nabla}f(x_{t-1})- \eta Q (x_t-x_{t+1}) + \eta Q (x_{t-1}-x_t)\|^2 \nonumber \\
 & = \sigma_{A}^{-1}\|\hat{\nabla}f(x_{t})-\nabla f(x_{t}) + \nabla f(x_{t}) - \nabla f(x_{t-1}) + \nabla f(x_{t-1})
    - \hat{\nabla}f(x_{t-1}) \nonumber \\
 & \quad -\eta Q (x_t-x_{t+1}) + \eta Q (x_{t-1}-x_t)\|^2 \nonumber \\
 & \mathop{\leq}^{i} \frac{5}{\sigma_{A}} \|\hat{\nabla}f(x_{t})-\nabla f(x_{t})\|^2
 + \frac{5}{\sigma_{A}} \|\hat{\nabla}f(x_{t-1})-\nabla f(x_{t-1})\|^2
 + \frac{5\eta^2\phi^2_{\max}}{\sigma_{A}} \|x_t-x_{t+1}\|^2 \nonumber \\
 & \quad + \frac{5(\eta^2\phi^2_{\max}+L^2)}{\sigma_{A}} \|x_{t-1}-x_{t}\|^2,
 \end{align}
 where the inequality (i) holds by the Assumption 1.

 Taking expectation conditioned on information $i_t$ to (78), we have
 \begin{align}
\mathbb{E}\|\lambda_{t+1}-\lambda_{t}\|^2 &\leq \frac{5}{\sigma_{A}} \mathbb{E} \|\hat{\nabla}f(x_{t})-\nabla f(x_{t})\|^2
  + \frac{5}{\sigma_{A}} \mathbb{E} \|\hat{\nabla}f(x_{t-1})-\nabla f(x_{t-1})\|^2 \nonumber \\
 & \quad + \frac{5\eta^2\phi^2_{\max}}{\sigma_{A}} \|x_t-x_{t+1}\|^2 + \frac{5(\eta^2\phi^2_{\max}+L^2)}{\sigma_{A}} \|x_{t-1}-x_{t}\|^2 \nonumber \\
 & \mathop{\leq}^{(i)} \frac{5L^2}{\sigma_{A} n} \sum_{i=1}^n\mathbb{E}\|x_{t}-z^{t}_i\|^2
 + \frac{5L^2}{\sigma_{A} n}\sum_{i=1}^n \mathbb{E} \|x_{t-1}-z^{t-1}_i\|^2
 + \frac{5\eta^2\phi^2_{\max}}{\sigma_{A}}\|x_{t+1}-x_t\|^2 \nonumber \\
 & \quad + \frac{5(\eta^2\phi^2_{\max} + L^2)}{\sigma_{A}}\|x_{t}-x_{t-1}\|^2, \nonumber
 \end{align}
 where the inequality $(i)$ holds by the Lemma 3.
\end{proof}

\section{Proof of the Lemma 14}
\label{app:lemma 14}

\begin{proof}
By the step 5 of Algorithm 4, we have
\begin{align}
 \mathcal {L}_\rho (x_t, y_{t+1},\lambda_t) \leq \mathcal {L}_\rho (x_t, y_{t},\lambda_t).
\end{align}
Next, by the optimal condition of step 7 in Algorithm 4, we have
\begin{align}
0 & =(x_t-x_{t+1})^T\big[\hat{\nabla}f(x_t)+\rho A^T(Ax_{t+1}+By_{t+1}-c) -A^T\lambda_t-\eta Q(x_t-x_{t+1})\big] \nonumber \\
 & = (x_t-x_{t+1})^T\big[\hat{\nabla}f(x_t) - \nabla f(x_t) + \nabla f(x_t) -A^T\lambda_t - \eta Q(x_t-x_{t+1})
 +\rho A^T(Ax_{t+1}+By_{t+1}-c)\big] \nonumber \\
 & \mathop{\leq}^{(i)} f(x_t) - f(x_{t+1}) + (x_t -x_{t+1})^T(\hat{\nabla}f(x_t) - \nabla f(x_t)) + \frac{L}{2}\|x_{t+1}-x_t\|^2 -\eta\|x_{t+1}-x_t\|^2_Q  \nonumber \\
 &\quad - \lambda_t^T(Ax_{t+1}-A x_t) + \rho (Ax_t -Ax_{t+1})^T(Ax_{t+1}+By_{t+1}-c) \nonumber \\
 & \mathop{=}^{(ii)} f(x_t) - f(x_{t+1}) + (x_t -x_{t+1})^T(\hat{\nabla}f(x_t) - \nabla f(x_t)) + \frac{L}{2}\|x_{t+1}-x_t\|^2 -\eta\|x_{t+1}-x_t\|^2_Q \nonumber \\
 & \quad - \lambda_t^T(A x_t+By_{t+1}-c) + \lambda_t^T(Ax_{t+1}+By_{t+1}-c)
  + \frac{\rho}{2}\|Ax_{t}+By_{t+1}-c\|^2 \nonumber \\
 & \quad - \frac{\rho}{2}\|Ax_{t+1}+By_{t+1}-c\|^2 -\frac{\rho}{2}\|Ax_t-Ax_{t+1}\|^2 \nonumber \\
 & = \mathcal {L}_\rho (x_t, y_{t+1},\lambda_t)- \mathcal {L}_\rho (x_{t+1}, y_{t+1},\lambda_t)
 + (x_t -x_{t+1})^T(\hat{\nabla}f(x_t) - \nabla f(x_t)) \nonumber \\
 & \quad + \frac{L}{2}\|x_{t+1}-x_t\|^2 -\eta\|x_{t+1}-x_t\|^2_Q -\frac{\rho}{2}\|Ax_t-Ax_{t+1}\|^2 \nonumber \\
 & \mathop{\leq}^{(iii)} \mathcal {L}_\rho (x_t, y_{t+1},\lambda_t)- \mathcal {L}_\rho (x_{t+1}, y_{t+1},\lambda_t)
 + (x_t -x_{t+1})^T(\hat{\nabla}f(x_t) - \nabla f(x_t)) \nonumber \\
 & \quad - (\eta\phi_{\min} + \frac{\sigma_A\rho}{2}-\frac{L}{2} )\|x_t-x_{t+1}\|^2,
\end{align}
where the inequality (i) holds by the Assumption 1; the equality (ii) holds by using the equality
$(a-b)^T(b-c) = \frac{1}{2}(\|a-c\|^2-\|a-b\|^2-\|b-c\|^2)$ on the term $\rho (Ax_t -Ax_{t+1})^T(Ax_{t+1}+By_{t+1}-c)$;
the inequality (iii) holds by using $-\phi_{\min}\|x_{t+1}-x_t\|^2 \geq -\|x_{t+1}-x_t\|^2_Q$ and $-\sigma_A\|x_{t+1}-x_t\|^2 \geq -\|Ax_t -Ax_{t+1}\|^2$.
Taking expectation conditioned on information $i_t$ to (80), and using $\mathbb{E}[\hat{\nabla} f(x_{t})]=\nabla f(x_{t})$, we have
\begin{align}
\mathbb{E} [\mathcal {L}_\rho (x_{t+1}, y_{t+1},\lambda_t)] \leq \mathcal {L}_\rho (x_t, y_{t+1},\lambda_t) - (\eta\phi_{\min}
+ \frac{\sigma_A\rho}{2} -\frac{L}{2})\|x_t-x_{t+1}\|^2.
\end{align}

By the step 8 of Algorithm 4, and taking expectation conditioned on information $i_t$, we have
\begin{align}
\mathbb{E} [\mathcal {L}_\rho (x_{t+1}, y_{t+1},\lambda_{t+1})-\mathcal {L}_\rho (x_{t+1}, y_{t+1},\lambda_t)]
&= \frac{1}{\rho}\mathbb{E} \|\lambda_t-\lambda_{t+1}\|^2 \nonumber \\
& \mathop{\leq}^{(i)} \frac{5L^2}{\sigma_{A} \rho n} \sum_{i=1}^n \mathbb{E}\|x_{t}-z^{t}_i\|^2
 + \frac{5L^2}{\sigma_{A}\rho n}\sum_{i=1}^n \mathbb{E}\|x_{t-1}-z^{t-1}_i\|^2\nonumber \\
& \quad + \frac{5(\eta^2\phi^2_{\max})}{\sigma_{A} \rho}\|x_{t+1}-x_t\|^2 + \frac{5(\eta^2\phi^2_{\max}+L^2)}{\sigma_{A} \rho}\|x_{t}-x_{t-1}\|^2,
\end{align}
where the inequality $(i)$ holds by the Lemma 13.

Combining (79), (81) and (82), we have
\begin{align}
 \mathbb{E} [\mathcal {L}_\rho (x_{t+1}, y_{t+1},\lambda_{t+1})] \leq  & \mathcal {L}_\rho (x_t, y_{t},\lambda_t) + \frac{5L^2}{\sigma_{A}\rho n} \sum_{i=1}^n \mathbb{E}\|x_{t}-z^{t}_i\|^2
  + \frac{5L^2}{\sigma_{A}\rho n}\sum_{i=1}^n \mathbb{E}\|x_{t-1}-z^{t-1}_i\|^2  \nonumber \\
 & + \frac{5(\eta^2\phi^2_{\max}+L^2)}{\sigma_{A}\rho}\|x_{t}-x_{t-1}\|^2
  - \big(\eta\phi_{\min} + \frac{\sigma_A\rho}{2} - \frac{L}{2} - \frac{5\eta^2\phi^2_{\max}}{\sigma_{A}\rho} \big)\|x_{t+1}-x_t\|^2.
\end{align}
Next, considering $\frac{1}{n}\sum_{i=1}^n \mathbb{E}\|x_{t+1}-z^{t+1}_i\|^2$, we have
\begin{align}
\frac{1}{n}\sum_{i=1}^n \mathbb{E}\|x_{t+1}-z^{t+1}_i\|^2 =\frac{1}{n}\sum_{i=1}^n \big[ \frac{1}{n}\mathbb{E}\|x_{t+1}-x_{t}\|^2 + \frac{n-1}{n}\mathbb{E}\|x_{t+1} - z^{t}_i\|^2 \big].
\end{align}
The term $\mathbb{E}\|x_{t+1} - z^{t}_i\|^2$ in (84) can be bounded as follows:
\begin{align}
\mathbb{E}\|x_{t+1} - z^{t}_i\|^2 & = \mathbb{E}\|x_{t+1}-x_{t}+x_{t}-z^{t}_i\|^2 \nonumber \\
& = \mathbb{E}[\|x_{t+1}-x_{t}\|^2+ 2(x_{t+1}-x_{t})^T(x_{t}-z^{t}_i ) +\|x_{t-1}-z^{t}_i\|^2] \nonumber \\
& \leq \mathbb{E}[\|x_{t+1}-x_{t}\|^2 + 2(\frac{1}{2\beta}\mathbb{E}\|x_{t+1}-x_{t}\|^2 + \frac{\beta }{2}\|x_{t}-z^{t}_i\|^2 )
 +\|x_{t}-z^{t}_i\|^2] \nonumber \\
& = (1+\frac{1}{\beta})\mathbb{E}\|x_{t+1}-x_{t}\|^2 + (1+\beta)\|x_{t}-z^{t}_i\|^2, \nonumber
\end{align}
where $\beta>0$, and the inequality is due to the Cauchy inequality.
Then, we have
\begin{align}
 &\frac{1}{n}\sum_{i=1}^n \mathbb{E}\|x_{t+1}-z^{t+1}_i\|^2 \leq (1+\frac{(n-1)}{n\beta})\mathbb{E}\|x_{t+1}-x_{t}\|^2  +(1+\beta)\frac{n-1}{n^2}\sum_{i=1}^n \|x_{t}-z^{t}_i\|
\end{align}
Combining the inequalities (83) and (85), we have
\begin{align}
 &\mathbb{E} \big[\mathcal {L}_\rho (x_{t+1}, y_{t+1},\lambda_{t+1}) + \frac{5(L^2 + \eta^2\phi^2_{\max})}{\sigma_{A} \rho}\|x_{t+1}-x_t\|^2
  + \frac{\alpha_{t+1}}{n}\sum_{i=1}^n (\|x_{t+1}-z^{t+1}_i\|^2 + \|x_{t}-z^{t}_i\|^2) \big ] \nonumber \\
 & \leq \mathcal {L}_\rho (x_{t}, y_{t},\lambda_{t}) + \frac{5(L^2+\eta^2\phi^2_{\max})}{\sigma_{A} \rho} \|x_t - x_{t-1}\|^2 \nonumber\\
 & \quad + \big[(2+\beta-\frac{1+\beta}{n})\alpha_{t+1} + \frac{5L^2}{\sigma_{A} \rho}\big]\frac{1}{n}\sum_{i=1}^n \big( \|x_{t}-z^{t}_i\|^2
  + \|x_{t-1}-z^{t-1}_i\|^2 \big) \nonumber\\
 & \quad -\big[\eta\phi_{\min} + \frac{\sigma_A\rho}{2}- \frac{L}{2}-\frac{5(2\eta^2\phi^2_{\max} + L^2)}{\sigma_{A} \rho}
   -(1+\frac{1}{\beta}-\frac{1}{n\beta})\alpha_{t+1}\big]\|x_{t+1}-x_t\|^2 \nonumber\\
 & \quad -(2+\beta-\frac{1+\beta}{n})\frac{\alpha_{t+1}}{n}\sum_{i=1}^n \|x_{t-1}-z^{t-1}_i\|^2.
\end{align}

Finally, using the definition of the sequence $\{\hat{\Phi}\}_{t=1}^T$, (28) and (29), we have
\begin{align}
\hat{\Phi}_{t+1}  \leq  \hat{\Phi}_{t}-\Gamma_{t}\|x_{t+1}-x_t\|^2
  -(2+\beta-\frac{1+\beta}{n})\frac{\alpha_{t+1}}{n}\sum_{i=1}^n \|x_{t-1}-z^{t-1}_i\|^2.
\end{align}
Since $\Gamma_t >0$, we can obtain the above result.
\end{proof}

%
%

\section{Proof of the Theorem 16}
\label{app:theorem 16}

\begin{proof}
Using the above inequality (87), we have
\begin{align}
 \hat{\Phi}_{t+1} & \leq \hat{\Phi}_{t} -\Gamma_{t+1}\|x_{t+1}-x_t\|^2
  - (2+\beta-\frac{1+\beta}{n})\frac{\alpha_{t+1}}{n}\sum_{i=1}^n\|x_{t-1}-z^{t-1}_i\|^2,
\end{align}
for $t\in\{1,2,\cdots,T\}$.
Summing the inequality (88) over $t=1,2,\cdots,T$, then we have
\begin{align}
 \hat{\Phi}_{T}  \leq \hat{\Phi}_{1} - \gamma\sum_{t=1}^{T} \mathbb{E}\|x_{t+1}-x_t\|^2
   - \omega \sum_{t=1}^T\frac{1}{n}\sum_{i=1}^n\|x_{t-1}-z^{t-1}_i\|^2.
\end{align}
where $\gamma = \min_t \Gamma_t $ and $\omega=\min_t (2+\beta-\frac{1+\beta}{n})\alpha_{t+1}$.
From Lemma 15, there exists a lower bound $\hat{\Phi}^*$ of the sequence $\{\hat{\Phi}_t\}$, i.e.,
$ \hat{\Phi}_{t} \geq \hat{\Phi}^*$ holds for $\forall t \geq 1$.
By the definition of $\hat{\theta}_{t}$ and (89), we have
\begin{align}
\hat{\theta}_{\tilde{t}}=\min_{1 \leq t \leq T} \hat{\theta}_{t} \leq \frac{1}{\tau T} (\hat{\Phi}_{1} - \hat{\Phi}^*),
\end{align}
where $\tau=\min(\gamma,\omega)$, so $\hat{\theta}_{\tilde{t}} = O(\frac{1}{T})$.

Next, we give upper bounds to the terms in (11-13) by using $\hat{\theta}_{t}$.
Using (77), we have
\begin{align}
 & \mathbb{E}\|A^T\lambda_{t+1}-\nabla f(x_{t+1})\|^2  \nonumber \\
 & = \mathbb{E}\|\hat{\nabla}f(x_{t}) - \nabla f(x_{t+1}) - \eta Q (x_t-x_{t+1})\|^2 \nonumber \\
 & = \mathbb{E}\|\hat{\nabla}f(x_{t})-\nabla f(x_{t}) +\nabla f(x_{t})- \nabla f(x_{t+1})
  - \eta Q (x_t-x_{t+1})\|^2  \nonumber \\
 & \leq\frac{3L^2}{n}\sum_{i=1}^n\|x_{t}-z^t_i\|^2  + 3(L^2+\eta^2\phi^2_{\max})\|x_t-x_{t+1}\|^2 \nonumber \\
 & \leq 3(L^2+\eta^2\phi_{\max}^2)\theta_{t}.
\end{align}
By the step 8 of Algorithm 4 and the Lemma 13, we have
\begin{align}
\mathbb{E}\|Ax_{t+1}+By_{t+1}-c\|^2 &= \frac{1}{\rho^2}\|\lambda_{t+1}-\lambda_t\|^2  \nonumber \\
 & \leq \frac{5L^2}{\sigma_{A}\rho^2 n} \sum_{i=1}^n \mathbb{E} \|x_{t+1}-z^{t+1}_i\|^2
 + \frac{5L^2}{\sigma_{A}\rho^2 n} \sum_{i=1}^n \|x_{t}-z^t_i\|^2  \nonumber \\
 & \quad + \frac{5L^2}{\sigma_{A}\rho^2} \mathbb{E} \|x_{t+1}-x_t\|^2
  + \frac{5(L^2+\eta^2\phi^2_{\max})}{\sigma_{A}\rho^2}\|x_{t}-x_{t-1}\|^2  \nonumber \\
 & \leq \frac{5(L^2+\eta^2\phi_{\max}^2)}{\sigma_{A}\rho^2} \theta_{t}.
\end{align}
By the step 5 of Algorithm 4,
there exists a subgradient $\mu \in \partial g(y_{t+1})$ such that
\begin{align}
 \mathbb{E}\big[\mbox{dist} (B^T\lambda_{t+1}, \partial g(y_{t+1}))\big]^2  & \leq \|\mu-B^T\lambda_{t+1}\|^2 \nonumber \\
 & = \|B^T\lambda_t-\rho B^T(Ax_t+By_{t+1}-c)-B^T\lambda_{t+1}\|^2 \nonumber \\
 & = \|\rho B^TA(x_{t+1}-x_{t})\|^2 \nonumber \\
 & \leq \rho^2\|B\|_2^2\|A\|_2^2\|x_{t+1}-x_t\|^2 \nonumber \\
 & \leq \rho^2\|B\|_2^2\|A\|_2^2 \theta_{t}.
\end{align}
Thus, by (32) and the Definition 4, we conclude that the SAGA-ADMM can converge an $\epsilon$-stationary point of
the problem (1).
\end{proof}

\section{Proof of the Theorem 17}
\label{app:theorem 17}

\begin{proof}
 By the optimal condition of the step 6 in Algorithm 1, we have
 \begin{align}
  0 & = \nabla f_{i_t}(x_t)-A^T\lambda_t+\rho A^T(Ax_{t+1}+By_{t+1}-c) - \eta Q(x_t - x_{t+1}) \nonumber \\
    & = \nabla f_{i_t}(x_t)- A^T\lambda_{t+1} - \eta Q (x_t - x_{t+1}),
 \end{align}
 where the second equality holds by the step 7 in Algorithm 1.
 By (94), we have
 \begin{align}
 \mathbb{E}\|A^T\lambda_{t+1}-\nabla f(x_{t+1})\|
 & = \mathbb{E}\|\nabla f_{i_t}(x_t) - \nabla f(x_{t+1}) - \eta Q (x_{t}-x_{t+1})\| \nonumber \\
 & = \mathbb{E}\|\nabla f_{i_t}(x_t) - \nabla f(x_{t}) + \nabla f(x_{t}) - \nabla f(x_{t+1})
  - \eta Q (x_{t}-x_{t+1})\| \nonumber \\
 & \geq \mathbb{E}\|\nabla f_{i_t}(x_t) - \nabla f(x_{t})\| - \mathbb{E}\|\nabla f(x_{t}) - \nabla f(x_{t+1})
  -\eta Q(x_{t}-x_{t+1})\| \nonumber \\
 & \mathop{\geq}^{(i)} \mathbb{E}\|\nabla f_{i_t}(x_t) - \nabla f(x_{t})\| - (L + \eta\phi_{\max})\|x_{t}-x_{t+1}\|  \nonumber \\
 & \geq \delta  - (L + \eta\phi_{\max})\|x_{t}-x_{t+1}\|
 \end{align}
 where the equality $(i)$ holds by the Assumption 1.
 Suppose the sequence $\{x_t,y_t,\lambda_t\}$ generated by Algorithm 1, and it
 converges to a stationary point $(x^*,y^*,\lambda^*)$
 of the problem (1).
 Then $\exists \epsilon >0$, we have
 \begin{align}
 \|x_{t+1}-x_t\| \leq \|x_{t+1}-x^*\| + \|x_{t}-x^*\| \leq \epsilon/2 + \epsilon/2 = \epsilon. \nonumber
 \end{align}
 If $\delta \geq 2(L + \eta\phi_{\max})\epsilon$, we have
 $\mathbb{E}\|A^T\lambda_{t+1}-\nabla f(x_{t+1})\| \geq (L + \eta\phi_{\max})\epsilon $.
 Thus, we obtain the above conclusion by contradiction.
\end{proof}




%


\end{document}